\documentclass[11pt,a4paper,openany]{article}

\usepackage{amsmath,amsfonts,ams}
\usepackage{latexsym}
\usepackage{amssymb}
\usepackage{color}

\numberwithin{equation}{section}
\newtheorem{proposition}{Proposition}[section]
\newtheorem{definition}{Definition}[section]
\newtheorem{lemma}{Lemma}[section]
\newtheorem{theorem}{Theorem}[section]

\newtheorem{Rem}{Remark}[section]

\begin{document}
\title{On Global Solution to the Klein-Gordon-Hartree Equation below Energy Space}
\author{Changxing Miao$^1$ and  Junyong Zhang$^2$\\
{\small  $^{1}$Institute of Applied Physics and Computational Mathematics,}\\
    {\small P.O. Box 8009, Beijing 100088, P.R. China.}\\
     E-mail:\quad {\small miao\_changxing@iapcm.ac.cn}\\
{\small $^{2}$The Graduate School of China Academy of Engineering Physics}\\
        {\small P. O. Box 2101,\;  Beijing,\; China,\; 100088 } \\
E-mail: {zhangjunyong111@sohu.com }
 \date{}
 }
\maketitle

\begin{abstract}
In this paper, we consider the Cauchy problem for Klein-Gordon
equation with a cubic convolution nonlinearity in $\R^3$. By making
use of Bourgain's method in conjunction with a precise Strichartz
estimate of S.Klainerman and D.Tataru, we establish the $H^s (s<1)$
global well-posedness of the Cauchy problem for the cubic
convolution defocusing Klein-Gordon-Hartree equation. Before
arriving at the previously discussed conclusion, we obtain global
solution for this non-scaling equation with small initial data in
$H^{s_0}\times H^{s_0-1}$ where $s_0=\frac\gamma 6$ but not
$\frac\gamma2-1$, for this equation that we consider is a
subconformal equation in some sense. In doing so a number of
nonlinear prior estimates are already established by using Bony's
decomposition, flexibility of Klein-Gordon admissible pairs which
are slightly different from that of wave equation and a commutator
estimate. We establish this commutator estimate by exploiting
cancellation property and utilizing Coifman and Meyer multilinear
multiplier theorem. As far as we know, it seems that this is the
first result on low regularity for this Klein-Gordon-Hartree
equation.
\end{abstract}

\begin{center}
 \begin{minipage}{120mm}
 \end{minipage}
 \end{center}

\section{Introduction}
We study the following Cauchy problem for the Klein-Gordon-Hartree
equation:
\begin{equation}\label{1}
\begin{cases}
\square \phi+\phi+(|x|^{-\gamma}\ast|\phi|^2)\phi=0~~~ \mathrm{in}~~~\R\times\R^3\\
\phi|_{t=0}=\phi_0,\quad\partial_t \phi|_{t=0}=\phi_1.\end{cases}
\end{equation}
Here $\phi(t,x)$ is a complex valued function defined in space time
$\R^{1+3}$, and $\square=\partial_{tt}-\Delta$.

Recently the Cauchy problem \eqref{1} has been extensively studied
in the case with initial data $(\phi_0,\phi_1)\in H^1(\R^n)\times
L^2(\R^n)$. The well-posedness and the asymptotic behavior of
solution to the Cauchy problem \eqref{1} have been studied by G.P.
Menzala and W.Strauss \cite{M82S, M83}. The scattering theory of
solution to \eqref{1} has been established in \cite{Mi07}. On the
other hand, the time-dependent Schr\"{o}dinger equation with
interaction term $(|x|^{-\gamma}\ast|\phi|^2)\phi$  has also been
extensively studied. Ginibre and Velo \cite{GiV00} gave the
scattering theory of Hartree equation for the energy subcritical
case. For the energy critical case and mass critical, one can refer
to \cite{MXZ07,MXZ09} with radial initial data.

Many authors \cite{B099, G03P, K00P, M02Z, T99} have studied the
local well-posedness (as well as global well-posedness) in
fractional Sobolev spaces for the Cauchy problem of general
semilinear wave or Schr\"odinger equations under minimal
regularity assumptions on the initial data. For example, Tao
\cite{T99} established the sharp local well-posedness of nonlinear
wave equation. Kenig, Ponce, and Vega \cite{K00P} had established
the global well-posedness under the energy norm for the Cauchy
problem of nonlinear wave equations with rough initial data (in
particular, in $\dot H^s(\R^3),~\frac34<s<1$ for cubic wave
equation). They used the Fourier truncation method discovered by
Bourgain \cite{B099}. And also \cite{M02Z} extended
Kenig-Ponce-Vega's result to the dimension $n\geq4$. Recently, I.
Gallagher and F. Planchon \cite{G03P} presented a different proof
of the result in \cite{K00P} for $\frac34 < s < 1$. H. Bahouri and
Jean-Yves Chemin \cite{B06C} proved global well posedness for $s =
\frac34$ by using a nonlinear interpolation method and logarithmic
estimates from S. Klainermann and D. Tataru\cite{K99T}. We also
find Roy \cite{Roy07} obtains the global well-posedness for rough
initial data in $\dot H^s, \frac23<s<1$ by following the
$I$-method \cite{C02T} and scaling transformation. However, if one
similarly deals with Klein-Gordon equation by using $I$-method, he
or she may meet a problem caused by the lack of the scaling
property. More studies and discussions on the low regularity of
nonlinear wave or dispersive Schr\"odinger equations could be
found in \cite{B099, T06}. However, as far as we know, very few
authors are engaged in studying the global well-posedness of the
Cauchy problems \eqref{1} with less regular initial data. It is
natural to ask whether a similar or better result holds for the
problem \eqref{1}.

This paper endeavors to find a global well-posedness solution to the
Cauchy problem \eqref{1} with initial data $(\phi_0,\phi_1)\in
H^s(\R^3)\times H^{s-1}(\R^3)$ for some $s>\frac\gamma4$ with
$\gamma\in(2,3)$. Now we should remark some differences between
\eqref{1} and cubic wave equation. If one views \eqref{1} as a wave
equation by dropping the massive term and then makes some scaling
analysis, we will find this nonlocal nonlinear term shares the
scaling property of the nonlinearity $|u|^{\frac4{5-\gamma}}u$. One
can check that $k:=\frac4{5-\gamma}+1<3$ when $2<\gamma<\min\{n,4\}$
with $n=3$ and this result shows that the equation which we consider
is in subconformal case. To obtain the global well-posedness theory,
some previous literatures also show the subconformal equations are
slightly different from the superconformal ones. For instance,
Lindblad and Sogge \cite{L95S} \cite{S95} have shown the global
existence and scattering theory for small data in a less regularity
space for the superconformal case, while not for the subconformal
case. Inspired by \cite{G03P}, we also split the initial data into
low frequency part data in $H^1$ and high frequency part data in
$H^{s_0}$ with a suitable $s_0$. Since the problem \eqref{1} is
global well-posed for large data in $H^1$ and small data in
$H^{s_0}$, one may be tempted to follow a general principle of
nonlinear interpolation and claim the problem \eqref{1} is global
well-posed between them. Compared with the cubic wave equation,
speaking of the Strichartz estimate, we believe that the global
solution with high frequency data should exist in
$H^{\frac\gamma2-1}$. It is well known that the Strichartz estimate
is associated with scaling transform and it is scaling invariant.
Unfortunately, the equation that we consider is a subconformal one,
and its concentration effects take over scaling. Since the
Strichartz estimate is applied to our subconformal equation, hence
this brings about some loss to get a better result. In order to get
a better result, one should establish an estimate which is conformal
invariant. Fortunately, we can take $ 0\leq\theta\leq 1$ as a
parameter for the flexible admissible pairs (see Definition
\ref{def2.5})to make the Strichartz estimate of Klein-Gordon more
flexible than wave equation. This helps us to get a global solution
with the high frequency data, at the cost of
$0\leq\theta=\frac{6}\gamma-2\leq1$ which weakens the Strichartz
estimate and causes $2<\gamma<3$. One can refer the detail in
Section 3.

We point out that it is easy to have the result for
$\frac{\gamma}3-\frac16<s<1$ by rough H\"{o}lder's inequality. But
how to get our low bound $ \frac\gamma4<s<1$? A good way to think
about this is via precise Strichartz estimate to obtain index $s$ as
low as possible. The nonlinearity including a formal negative
derivative brings us some difficulties caused by the fact that the
negative derivative acts on the low frequency part. And this leads
us to restricts $s>\frac\gamma4$ rather than $s>\max\{\frac12,
\frac\gamma2-\frac34\}$. At the end of this section, we also give
some intuitive analysis to show our result is reasonable. As a
limited case, our result recovers the result of \cite{G03P, K00P}
when $\gamma$ tends to $3$.

During the process of proving our key estimate Lemma \ref{lem6}, the
nonlocal nonlinearity brings about some essential difficulties when
we try to make use of the precise Strichartz estimate. Compared with
the general semilinear nonlinearity, the convolution nonlinearity
not only essentially represents a negative derivation in it but also
has a difference construction of nonlinearity. These differences and
difficulties prevent us from obtaining directly our expected result
$s>\frac\gamma4$ by restricting the range of the parameter $r$. To
overcome these difficulties, we firstly construct a commutator and
establish this commutator estimate by exploiting cancellation
property and utilizing Coifman and Meyer multilinear multiplier
theorem and then go on our process through using precise Strichartz
estimate.

Now we state our main result:
\begin{theorem}\label{the1}
Let $\frac\gamma 4<s<1$ with $2<\gamma<3$. If $(\phi_0,\phi_1)\in
H^s(\R^3)\times H^{s-1}(\R^3)$, then there exists  a unique global
 solution $\phi$ of \eqref{1} such that
$\phi\in C(\R^+;H^s(\R^3))$.
\end{theorem}

We conclude this section by giving a sketch of the proof of Theorem
\ref{the1} and one shall read more detailed information in the rest
of this paper. Without loss of generality, we only consider $\phi$
as a real function for simplicity from now on. Since the problem
\eqref{1} is global well-posed for large data in $H^1$ and small
data in $H^{s_0}$ with $s_0=\frac{\gamma}{6}$, one may be tempted to
follow a general principle of nonlinear interpolation and believe
the problem \eqref{1} to be global well-posed between them, as well
as the cubic defocusing wave equation \cite{G03P}. To make sense of
this heuristic, we proceed it in the following steps. \vskip0.3cm

Step 1.\;  The purpose of this step is to show the global
well-posedness for the high frequency part. We split the initial
data:
$$\phi_i=(\mathrm{I}-S_J)\phi_i+S_J\phi_i\stackrel{\mathrm{def}}{=}v_i+u_i\quad \quad i=0,1$$
where $\mathrm{I}$ is identity operator and $S_J$ is
Littlewood-Paley operator, referring  to Section 2. It is easy to
see that
$$\|u_0\|_{H^{1}}\lesssim2^{J(1-s)}\|\phi_0\|_{H^{s}},\quad\|u_1\|_{L^2}\lesssim\|\phi_1\|_{L^2}$$
and
$$\|v_0\|_{H^{\beta}}\lesssim2^{J(\beta-s)}\|\phi_0\|_{H^{s}},
\quad\|v_1\|_{H^{\beta-1}}\lesssim2^{J(\beta-s)}\|\phi_1\|_{H^{s-1}}~~~~~\mbox{for
all}~\beta\leq s.$$ Thus it follows that
\begin{align}
\label{3}\mathcal{E}_{h,\sigma}
\lesssim& \; 2^{J(\sigma-s)}\mathcal{E}_s,\quad\mbox{for}~~\sigma\leq s\\
\label{4}\mathcal{E}_{\ell,1}\lesssim&\;
2^{J(1-s)}\mathcal{E}_s,\quad\mbox{for}~~s\leq1,
\end{align}
where
\begin{align}
\label{5}\mathcal{E}_s\stackrel{\mathrm{def}}{=\!=}\|\phi_0\|_{H^{s}}+\|\phi_1\|_{H^{s-1}},\\
\label{6}\mathcal{E}_{h,\sigma}\stackrel{\mathrm{def}}{=\!=}\|v_0\|_{H^{\sigma}}+\|v_1\|_{H^{\sigma-1}},\\
\label{7}\mathcal{E}_{\ell,\sigma}\stackrel{\mathrm{def}}{=\!=}\|u_0\|_{H^{\sigma}}+\|u_1\|_{H^{\sigma-1}}.
\end{align}
Choosing $J$ large enough, one can achieve $\mathcal{E}_{h,s_0}$
small enough, in other words, initial data of the following problem
\begin{equation}\label{8}
\begin{cases}
\square v+v+(|x|^{-\gamma}\ast v^2)v=0~~~~~ in~~~\R\times\R^3,\\
v|_{t=0}=v_0, \quad\partial_t v|_{t=0}=v_1\end{cases}
\end{equation}
is small enough in $H^{s_0}(\R^3)\times H^{s_0-1}(\R^3)$ where
$s_0<s$. Due to some technique difficulties and this equation is
subconformal one, we are restricted to choose $s_0=\frac\gamma6$
while not $\frac\gamma2-1$ proposed by scaling analysis or
$\frac\gamma4-\frac14$ proposed by conformal analysis. We will get a
global well-posed solution to the Cauchy problem \eqref{8}, see
Section 3 for details.

\vskip0.3cm
 Step 2. In order to recover a
solution to our problem \eqref{1}, we solve a perturbed equation
with large initial data in $H^1\times L^2$,
\begin{equation}\label{9}
\begin{cases}
\square u+u+\mathcal{I}(u^2)u+2\mathcal{I}(uv)u+\mathcal{I}(v^2)u
+\mathcal{I}(u^2)v+2\mathcal{I}(uv)v=0,\\
u|_{t=0}=u_{0}\quad\partial_t u|_{t=0}=u_{1},
\end{cases}
\end{equation}
where the operator $\mathcal{I}$ is the operator
$(-\Delta)^{\frac{\gamma-3}2}$. We will prove there exists a unique
local solution to \eqref{9} in $C([0,T];H^1)$. \vskip0.3cm

Step 3.\; To complete the proof of Theorem \ref{the1}, the key
 is how to  extend  the local solution to a
global solution. We should establish a priori bound on the energy of
the local solution $u$. In fact, the energy estimate yields
\begin{align*}
&\frac12\left(\|u(t)\|^2_{H^1}
+\|u_t(t)\|^2_{L^2}\right)+\frac12\int_{\R^3\times\R^3}|x-y|^{-\gamma}u^2(y,t)u^2(x,t)dydx\\
\leq &\frac12\left(\|u_0\|^2_{H^1}
+\|u_1\|^2_{L^2}\right)+\frac12\int_{\R^3\times\R^3}|x-y|^{-\gamma}u_0^2(y)u_0^2(x)dydx\\
&+\left|\int_{0}^t\int_{\R^3}
\mathcal{I}(v^2)(x,\tau)u(x,\tau)\partial_\tau u(x,\tau)dx
d\tau\right|\\
&+2\left|\int_{0}^t\int_{\R^3}
\mathcal{I}(uv)(x,\tau)v(x,\tau)\partial_\tau u(x,\tau)dx
d\tau\right|\\
&+\left|\int_{0}^t\int_{\R^3}
\mathcal{I}(u^2)(x,\tau)v(x,\tau)\partial_\tau u(x,\tau)dx
d\tau\right|\\
&+2\left|\int_{0}^t\int_{\R^3}
\mathcal{I}(uv)(x,\tau)u(x,\tau)\partial_\tau u(x,\tau)dx
d\tau\right|.
\end{align*}

Let $H_T(u):=\sup\limits_{t<T}H(u)(t)$ where
$$H(u)(t)\stackrel{\mathrm{def}}{=}\left(\frac12\|u(t)\|^2_{H^1}
+\frac12\|u_t(t)\|^2_{L^2}+\frac12\int_{\R^3\times\R^3}|x-y|^{-\gamma}u^2(y,t)u^2(x,t)dydx\right)$$
and then by making use of H\"{o}lder's inequality and Sobolev
embedding, it follows that
\begin{align*}
H_T(u)\lesssim & H(u)(0)+H_T(u)\int_0^T\left\|
v(\tau)\right\|^2_{L^{\frac6{4-\gamma}}}d\tau+H^{\frac32}_T(u)\int_0^T\left\|
v(\tau)\right\|_{L^{\frac6{7-2\gamma}}}d\tau\cr\lesssim &
H(u)(0)+H_T(u)T^{\frac{7-\gamma}6}\|v\|^2_{
X^{\beta}}+H^{\frac32}_T(u)T^{\frac{5-\gamma}3}\|v\|_{
X^{\alpha}}\cr\lesssim &
2^{2J(1-s)}+H_T(u)T^{\frac{7-\gamma}6}2^{2J(\beta-s)}+H^{\frac32}_T(u)T^{\frac{5-\gamma}3}2^{J(\alpha-s)}
\end{align*}
where $\alpha=\frac{2\gamma-4}3, \beta=\frac{\gamma-1}3$ and the
space $X^\alpha$ is defined in the coming section. What we want to
do is to control $H_T(u)$ for arbitrarily large $T$. As long as
$s>(\alpha+1)/2=\frac{\gamma}3-\frac16$, by choosing $J$ large
enough, bootstrap argument yields
$$H_T(u)\lesssim 2^{2J(1-s)}. $$
One can see that, if $s>\frac{\gamma}3-\frac16$, the argument is
trivial, since the above mentioned result can be deduced from some
rough estimates such as the H\"{o}lder estimate. On the other hand,
since the scaling suggests us that $X^{\frac\gamma2-1}$ is the
lowest regularity space which $v$ could belong to, it is tempting
and reasonable to believe that the best result obtained by this
method is $s>(\frac{\gamma}2-1+1)/2=\frac\gamma4$ instead of
$\alpha$ by $\frac\gamma2-1$. To obtain this optimal result
$s>\frac\gamma4$, we adopt some more sophisticated tools such as
precise Strichartz estimate, Bony's paraproduct estimates and twice
Bony's decomposition. This result is achieved under an assumption of
a core estimate which will be shown through the precise Strichartz
estimate and a commutator estimate.

The paper is organized as follows: In the coming section, we recall
some notations and recollect some well known results on Besov spaces
in conjunction with the Littlewood-Paley theory which will be used
in the course of the proofs. Meanwhile, we also introduce  the
precise Strichartz estimate. Section 3 provides the global
well-posedness of original equation evoking the high frequency part
of initial data in $H^{s_0}$. In Section 4, we prove prove a local
well-posedness of perturbed equation with the low frequency of the
initial data in $H^1$ by the standard fixed point theorem. In
Section 5, we give a energy estimate for the low frequency part
provided an assumption the key estimate in Lemma \ref{lem6}. We
extend the local well-posedness of the perturbed equation to
globally well posed by the bootstrap argument in Section 6. In the
final section, we prove our essential and key lemma by the precise
Strichartz estimate, commutator estimate and Coifman and Meyer
multiplier theorem.

\section{Preliminaries}
\vskip0.3cm In this section, we shall present some well-known facts
on the Littlewood-Paley theory and introduce some notations,
definitions and estimates which are needed in this paper.
 Let ${\cal S}(\R^3)$ be
the Schwarz class of rapidly decreasing functions. Given $f\in {\cal
S}(\R^3)$, its Fourier transform ${\cal F}f=\hat f$ is defined by
$$
\hat f(\xi)=(2\pi)^{-\frac{3}2}\int_{\R^3}e^{-ix\cdot
\xi}f(x)dx,\quad {\cal F}^{-1}f=\hat f(-\xi).
$$
Choose two nonnegative radial functions $\chi$, $\varphi \in {\cal
S}(\R^3)$ supported respectively in ${\cal B}=\{\xi\in\R^3,\,
|\xi|\le\frac{4}{3}\}$ and ${\cal C}=\{\xi\in\R^3,\,
\frac{3}{4}\le|\xi|\le\frac{8}{3}\}$ such that
$$\chi(\xi)+\sum_{j\ge0}\varphi(2^{-j}\xi)=1,\quad\xi\in\R^3,$$
$$\sum_{j\in\Z}\varphi(2^{-j}\xi)=1,\quad\xi\in\R^3\backslash \{0\},$$
and
$$ \text{\rm supp}~\varphi(2^{-j}\cdot)\cap
\text{\rm supp}~\varphi(2^{-j'}\cdot)=\emptyset,\quad |j-j'|\geq
2,$$
$$\text{\rm supp}~\chi(\cdot)\cap \text{\rm
supp}~\varphi(2^{-j}\cdot)=\emptyset,\quad j\geq 1.$$

Now we are in position to define  the the Littlewood-Paley operators
$S_j$, ${\dot S}_j$,   $\triangle_j$ and ${\dot\triangle}_j$
 which are used to define Besov
space.
\begin{equation*}\triangle_{j}u\stackrel{\mathrm{def}}{=\!=}
\left\{\begin{aligned}&\qquad 0,\qquad \qquad \qquad \qquad \qquad \qquad  j\leq -2,\\
&\mathcal{F}^{-1}\big(\chi(\xi)\hat u(\xi)\big),\qquad \qquad \qquad \qquad j=-1,\\
&2^{jn}\displaystyle\int_{\R^n}(\mathcal{F}^{-1}\varphi)(2^jy)u(x-y)dy,\quad
j\geq 0,
\end{aligned}\right.\end{equation*}
$$S_ju\stackrel{\mathrm{def}}{=\!=}\sum_{j^\prime\leq j-1}\triangle_{j^\prime}u=
2^{jn}\int_{\R^n}(\mathcal{F}^{-1}\chi)(2^jy)u(x-y)dy,$$
\begin{equation*}\dot\triangle_{j}u\stackrel{\mathrm{def}}{=\!=}
2^{jn}\int_{\R^n}(\mathcal{F}^{-1}\varphi)(2^jy)u(x-y)dy, \quad
j\in\Z,\qquad
\end{equation*}
$$\dot S_ju\stackrel{\mathrm{def}}{=\!=}\sum_{j^\prime\leq j-1}{\dot\triangle}_{j^\prime}u.\qquad\qquad\qquad\qquad
\qquad\qquad\qquad$$ One easily shows that ${\dot\triangle}_j= {\dot
S}_{j+1}-{\dot S_j}$ for $j\in\Z$ and
$${\triangle}_{-1}= S_{0},\qquad {\dot\triangle}_j= {\triangle}_j,\quad  j\ge 0.$$

Now we give the Littlewood-Paley's  description of the  Besov
spaces.
\begin{definition}\label{def2.1}Let
$s\in \R, 1\le p,q\le\infty$. The homogenous Besov space $\dot
{B}^s_{p,q}$ is defined by
$$\dot {B}^s_{p,q}=\{f\in {\cal Z}'(\R^3):   \|f\|_{\dot
{B}^s_{p,q}}<\infty\},$$
 where
$$\|f\|_{\dot{B}^s_{p,q}}=\left\{\begin{array}{l}
\displaystyle\bigg(\sum_{j\in\Z}2^{jsq}\|{\dot\triangle}_j
f\|_p^q\bigg)^{\frac 1
q},\quad \hbox{for}\quad q<\infty,\\
\displaystyle\sup_{j\in \Z}2^{js}\|{\dot\triangle}_jf\|_p, \quad
\hbox{ for} \quad q=\infty,
\end{array}\right.
$$
and ${\cal Z}'(\R^3)$ can be identified by the quotient space ${\cal
S}'/{\cal P}$ with the  space ${\cal P}$ of polynomials.
\end{definition}
\begin{definition}\label{def2.2}Let
$s\in \R, 1\le p,q\le\infty$. The inhomogeneous Besov space $
{B}^s_{p,q}$ is defined by
$${B}^s_{p,q}=\{f\in {\cal S}'(\R^3): \|f\|_{
{B}^s_{p,q}}<\infty\},$$
 where
$$\|f\|_{{B}^s_{p,q}}=\left\{\begin{array}{l}
\displaystyle\bigg(\sum_{j\ge 0}2^{jsq}\|{\triangle}_j
f\|_p^q\bigg)^{\frac 1
q}+\|S_0(f)\|_p,\quad \hbox{for}\quad q<\infty,\\
\displaystyle\sup_{j\ge 0}2^{js}\|{\triangle}_jf\|_p+\|S_0(f)\|_p,
\quad \hbox{ for} \quad q=\infty.
\end{array}\right.
$$
\end{definition}

\noindent If $s>0$, then ${B}^s_{p,q}=L^p\cap\dot{B}^s_{p,q}$ and
$\|f\|_{B^s_{p,q}}\approx\|f\|_{p}+\|f\|_{\dot{B}^s_{p,q}}.$ We
refer the reader to \cite{Ber,Chemin1, Mi04, Tr78} for  details.
\vskip0.2cm

In order to investigate the low regularity solution of the Cauchy
problem \eqref{1}, we require the use of the smoothing effect
described by the Strichartz estimates and precise Strichartz
estimates. For the purpose of conveniently making use of the
Strichartz estimate, we introduce the admissible definition and the
resolution space.
\begin{definition}\label{def2.5}  We shall say that a pair $(q,r)$ is
admissible, for $0\leq \theta\leq 1$, if
$$q,r\geq 2, \quad  \; (q,r,\theta)\neq(2,\infty,0)\quad and \quad \frac1 q+\frac{2+\theta}{2r}\leq\frac{2+\theta}4.$$
\end{definition}

\begin{Rem}\label{rem2.2} The above admissible pairs in Definition \ref{def2.5} is more flexible than wave
admissible pairs, since  $\theta$ can vary from $0$ to $1$.
Obviously, an admissible pair in  Definition \ref{def2.5}   will
become a wave admissible pair when $\theta=0$. When we consider the
global existence for the high frequency part, we shall use
$\theta=\frac6{\gamma}-2$ since the equation that we consider is a
subconformal one.
\end{Rem}

The resolution space is defined in the following way based on the
admissible definition.
\begin{equation*}
\begin{split} X^\mu(I):=\bigcap_{0\leq\theta\leq1}X_{\theta}^\mu(I)
\end{split}
\end{equation*}
where
\begin{equation*}
\begin{split} X_{\theta}^\mu(I):=\Big\{&u: u\in({C}\cap L^\infty)(I;
H^{\mu})\cap L^q(I; B^{\sigma}_{r,2})\\& \mbox{(q,r) is admissible},
\frac1q=(3+\theta)(\frac12-\frac1r)+\sigma-\mu\Big\}.
\end{split}
\end{equation*}

We go on this section by recalling the classical Strichartz estimate
and the precise Strichartz estimate. This kind of estimate goes back
to Strichartz \cite{St70}, and has been proved in its generality by
Ginibre and Velo \cite{G85V}, and Keel and Tao \cite{K98T}. The
Strichartz estimates for the Klein-Gordon equation by using the
above flexible admissible pairs can be found in \cite{MZF2004}.

\begin{proposition} \label{pro1}
Let $u$ be a solution of
$$\square u+u=f\quad\mathrm{in}\quad \R\times\R^3\quad \mathrm{with}\quad u|_{t=0}=u_0,~~\partial_t u|_{t=0}=u_1.$$
Then, for any admissible pairs $(q_1,r_1)$ and $(q_2,r_2)$, we have
that
\begin{align}\nonumber&\|\triangle_j
u\|_{L^{q_1}(L^{r_1})}+2^{-j}\|\partial_t \triangle_j
u\|_{L^{q_1}(L^{r_1})}\\
\leq& C 2^{j(\frac {3+\theta} 2 -\frac {3+\theta} {r_1}-\frac 1
{q_1})}(\|\triangle_j u_0\|_{L^2}+2^{-j}\|\triangle_j
u_1\|_{L^2})\nonumber\\
&+C2^{j[{(3+\theta)}(1-\frac 1 {r_1}-\frac 1 {r_2})-\frac 1
{q_1}-\frac 1
{q_2}-1]}\|\triangle_jf\|_{L^{q'_2}(L^{r'_2})}.\label{11}
\end{align}
\end{proposition}

We shall see that the classical Strichartz estimates are not enough
to control some nonlinearities, and this leads us to resort to the
following precise Strichartz estimates which were established by
S.Klainerman and D.Tataru\cite{K99T}.
\begin{proposition}\label{pro2}
Let $u$ be a solution of
$$\square u+u=0\quad with\quad u|_{t=0}=u_0,\; \; \partial_t u|_{t=0}=u_1.$$
Assume that the supports of the Fourier transform of $u_0$ and $u_1$
are included in a ball $B(\xi_j,h2^j)$ with
$|\xi_j|\in[2^{j-2},2^{j+2}]$ and $h<\frac 1 8$. Then we have that,
for any admissible couple $(q,r)$,
\begin{equation}\label{12}\|
u\|_{L^{q}(L^{r})}+2^{-j}\|\partial_t u\|_{L^{q}(L^{r})}\leq C
2^{j(\frac {3} 2 -\frac {3} {r}-\frac 1 {q})}h^{\frac12
-\frac1r}(\|u_0\|_{L^2}+2^{-j}\|u_1\|_{L^2}).
\end{equation}
\end{proposition}

Let us recall the Hardy-Littlewood-Sobolev inequality \cite{Mi04,
St93} and a proposition of contraction which is generalization of
Picard's theorem \cite{Chemin1}. We denote operator $\mathcal{I}$ by
\begin{align*}
\mathcal{I}u\stackrel
{\mathrm{def}}{=\!=}(-\Delta)^{\frac{\gamma-3}2}u=|x|^{-\gamma}\ast
u,
\end{align*}
then
\begin{align}\label{13}
\|\mathcal{I}u\|_{L^q(\R^3)}\leq C_{p,q}\|u\|_{L^p(\R^3)}
\end{align}
for
$$0<\gamma<3,\quad 1<p<q<\infty,\quad\mathrm{and}\quad\frac1q=\frac1p-\frac{3-\gamma}3.$$

\begin{proposition}\label{prop2.3}
Let $X$ be a Banach space and let $B:\ X\times X\times \cdots\times
X\rightarrow X$ be a $m$-linear continuous operator {\rm($m\ge2$)}
satisfying
\begin{equation*}
\|B(u_1,u_2,\cdots,u_m)\|_X\le M\|u_1\|_X \|u_2\|_X \cdots
\|u_m\|_X,\qquad \forall u_1,\ u_2,\ \cdots,\ u_m\in X
\end{equation*}
for some constant $M>0.$ Let $\varepsilon>0$ be such that
$m(2\varepsilon)^{m-1} M<1.$ Then for every $y\in X$ with
$\|y\|_X\le\varepsilon$ the equation
\begin{equation}\label{13add1}  u=y+B(u,u,\cdots,u)
\end{equation}
has a unique solution $u\in X$ satisfying that $\|u\|_X\le
2\varepsilon.$ Moreover, the solution $u$ continuously depends on
$y$ in the sense that, if $\|y_1\|_X\le\varepsilon$ and $v=
y_1+B(v,v,\cdots,v),$ $\|v\|_X\le2\varepsilon$ then
\begin{eqnarray}\label{13add2}
\|u-v\|_X\le\frac1{1-m(2\varepsilon)^{m-1}M}\|y-y_1\|_X.
\end{eqnarray}
\end{proposition}

For the sake of convenience, we conclude this section by giving some
notations. The solution $\phi$ to the Cauchy problem \eqref{1} is
given by the following integral equation:
\begin{eqnarray*}
\phi(t,x)=\dot{K}(t){\phi_0}+K(t){\phi_1}+B(\phi,\phi,\phi)\stackrel{\mathrm{def}}{=}\mathcal{T}\phi
\end{eqnarray*} where
$$K(t):=\frac{\sin(t\sqrt{I-\Delta})}{\sqrt{I-\Delta}},$$
$$B(u_1,u_2,u_3):=-\int_0^tK(t-\tau)(|x|^{-\gamma}\ast(u_1u_2))u_3d\tau.$$
Throughout this article we shall denote by the letter $C$ all
universal constant and $\varepsilon>0$ is a arbitrary small data. We
shall sometimes replace an inequality of the type $f\leq Cg$ by
$f\lesssim g$. Also, we shall denote by $(c_j)_{j\in\Z}$ any
sequence of norm less than $1$ in $\ell^2(\Z)$. \vskip0.3cm

\section{Global existence for the high frequency part}

\vskip0.3cm Let us consider the Cauchy problem with the high
frequency data,
\begin{equation}\label{14}
\begin{cases}
\square v+v+(|x|^{-\gamma}\ast v^2)v=0,~~~~~(t,x)\in\R\times\R^3\\
v|_{t=0}=v_0, \quad\partial_t v|_{t=0}=v_1, ~~~~~ ~~~x\in\R^3.
\end{cases}
\end{equation}
and then its integral formation becomes
\begin{align}\nonumber
v(t, x)=&{\dot K}(t)v_0(x)+K(t)v_1(x)-\int_0^tK(t-\tau)(|x|^{-\gamma}\ast v^2)v d\tau\\
\stackrel {\mathrm{def}}{=\!=}&{\dot K}(t)v_0(x)+K(t)v_1(x)+ B(v, v,
v). \label{14add1}\end{align} Our goal in this section  is to prove
the global well-posedness of \eqref{14} or \eqref{14add1}. More
precisely, we have the following proposition:
\begin{proposition}\label{prop3}
Let $s_0=\frac{\gamma}{6}$ and suppose that $(v_0,v_1)\in
H^\mu\times H^{\mu-1}$ for any $0\leq\mu\leq 1$. There exists a
constant $\varepsilon_0>0$ such that if
$$\|v_0\|_{ H^{s_0}}+\|v_1\|_{ H^{s_0-1}}\leq \varepsilon_0,$$
then there exists a unique global solution $v$ to \eqref{14} or
\eqref{14add1} in ${X}^{s_0}(\R)\cap{X}^{\mu}{(\R)}$. Moreover,
$$\|v\|_{X^\mu}\leq C_\mu\left(\|v_0\|_{H^{\mu}}
+\|v_1\|_{H^{\mu-1}}\right).$$
\end{proposition}
\begin{Rem}
We focus on $\mu=\frac{2\gamma-4}3$ and $\mu=\frac{\gamma-1}3$ in
the coming section.
\end{Rem}
It is well known that the global existence theory for small initial
data is a straightforward result of nonlinear estimate, thus how to
obtain a suitable nonlinear estimate is essential. Before proving
this proposition, we make some analysis on nonlinear estimate. As
mentioned in the introduction, the nonlocal nonlinearity shares the
scaling with a subconformal nonlinearity when $\gamma<3$ and this
may bring some troubles when we make a choice of a suitable
resolution space $X^{s_0}$. Take $ 0\leq\theta\leq 1$ as a parameter
in the flexible admissible pairs (see Definition \ref{def2.5}), and
we make analysis on the relationship between $\theta$ and $s_0$. The
Strichartz estimate, H\"{o}lder inequality and
Hardy-Littlewood-Sobolev inequality imply that, for $\sigma\leq0$,
\begin{equation*}
\begin{split}
\|B(v,v,v)\|_{X^{s_0}}\leq
\|(|x|^{-\gamma}\ast|v|^2)v\|_{L^{q_1'}(B^{-\sigma}_{r_1',2})}\leq
\|v\|_{L^{q_2}(B^{-\sigma}_{r_2,2})}\|v\|^2_{L^{q_3}(L^{r_3})},
\end{split}
\end{equation*}
with satisfying
\begin{equation*}
\begin{split}
\frac1{q_1}&=(3+\theta)(\frac12-\frac1{r_1})+\sigma+s_0-1\cr
\frac1{q_2}&=(3+\theta)(\frac12-\frac1{r_2})-\sigma-s_0\cr
\end{split}
\end{equation*}
and
\begin{equation*}
\begin{split}
\frac1{q_3}&=(3+\theta)(\frac12-\frac1{r_3})-s_0\cr
1&=\frac1{q_1}+\frac1{q_2}+\frac2{q_3}\cr
2&=\frac\gamma3+\frac1{r_1}+\frac1{r_2}+\frac2{r_3}\cr
\end{split}
\end{equation*}
then
\begin{equation*}
\begin{split}
s_0=\frac\gamma2-1+\frac{\gamma\theta}6.
\end{split}
\end{equation*}
We find the fact index $s_0$ is increasing when the parameter
$\theta$ increases. It is tempting to choose $\theta=0$ to get the
smallest $s_0=\frac\gamma2-1$ proposed by scaling. However, in
addition the admissible condition implies that
\begin{equation*}
\begin{split}
\frac2{q_1}\leq(2+\theta)(\frac12-\frac1{r_1})\cr
\frac2{q_2}\leq(2+\theta)(\frac12-\frac1{r_2})\cr
\frac2{q_3}\leq(2+\theta)(\frac12-\frac1{r_3})\cr
\end{split}
\end{equation*}
then a direction computation gives that
\begin{equation*}
\begin{split}
2(\frac1{q_1}+\frac1{q_2}+\frac2{q_3})\leq(2+\theta)(2-\frac1{r_1}-\frac1{r_2}-\frac2{r_3})\cr
\end{split}
\end{equation*}
which yields that
\begin{equation*}
\begin{split}
\frac3\gamma\leq1+\frac\theta2.
\end{split}
\end{equation*}
If we choose $\theta=0$, then we are forced to $\gamma\geq 3$ which
contradict with our requirement $\gamma<3$. But if we choose
$\theta=\frac6\gamma-2$ and then $s_0=\frac\gamma6$ and we are
allowed by $2\leq \gamma\leq 3$.

{\bf Proof of Proposition \ref{prop3}} Thanks to Strichartz
estimate, we have
\begin{equation*}
\begin{split}
\|B(v,v,v)\|_{X^\mu}\leq
\|(|x|^{-\gamma}\ast|v|^2)v\|_{L^{q_1'}(B^{-\sigma}_{r_1',2})}\leq
\|v\|_{L^{q_2}(B^{-\sigma}_{r_2,2})}\|v\|^2_{L^{q_3}(L^{r_3})},
\end{split}
\end{equation*}
where
\begin{equation*}
\begin{split}
\big(\frac1{q_1}, \frac1
{r_1}\big)=\bigg(\frac3{3+\gamma}(1-\mu-\sigma),
\frac12-\frac\gamma{3+\gamma}(1-\mu-\sigma)\bigg),
\end{split}
\end{equation*}
and
\begin{equation*}
\begin{split}
\big(\frac1{q_2}, \frac1
{r_2}\big)=\bigg(\frac{3(\mu+\sigma)}{3+\gamma},
\frac12-\frac{\gamma(\mu+\sigma)}{3+\gamma}\bigg),\quad
\big(\frac1{q_3}, \frac1
{r_3}\big)=\bigg(\frac{\gamma}{2(3+\gamma)},
\frac{9+3\gamma-\gamma^2}{6(3+\gamma)}\bigg).
\end{split}
\end{equation*}
When $0\leq\mu\leq \frac12+\frac\gamma6$, we choose $\sigma=0$;
while $\frac12+\frac\gamma6<\mu\leq 1$, we choose
$\sigma=\frac12+\frac\gamma6-\mu$.
 Thus,
\begin{equation}\label{7a}
\begin{split}
\|B(v,v,v)\|_{X^\mu}\leq \|v\|_{X^\mu}\|v\|^2_{X^{s_0}}.
\end{split}
\end{equation}
Combining this nonlinear estimate, the Proposition \ref{prop3}
follows from a standard contraction argument and small initial data
condition.

\section{Local existence for the low frequency part}

\vskip0.3cm In this part, we shall study the following perturbed
problem  in $\R\times\R^3$:
\begin{equation}\label{21}
\begin{cases}
\square u+u+\mathcal{I}(u^2)u+2\mathcal{I}(uv)u+\mathcal{I}(v^2)u
+\mathcal{I}(u^2)v+2\mathcal{I}(uv)v=0\\
u|_{t=0}=u_{0}\quad\partial_t u|_{t=0}=u_{1}.
\end{cases}
\end{equation}

\begin{proposition}\label{prop6}
\quad Let $\alpha=\frac{2\gamma-4}3, \beta=\frac{\gamma-1}3$ and
 assume that $v$ be in $
{X}^{\alpha}\cap {X}^{\beta}$ and $(u_{0},u_{1})\in H^1\times L^2$,
then there exists a positive time $T$ such that a unique solution
$u$ to \eqref{21} satisfying $$u\in C([0,T]; H^1).$$
\end{proposition}
{\bf Proof of the Proposition \ref{prop6}}\quad  In practice,
solving \eqref{21} on $[0,T]$ is equivalent to solving the following
integral equation
\begin{align*}
u=&\dot{K}(t)u_{0}+K(t)u_{1}\\
&+\int_0^t K(t-\tau)\Big[\mathcal{I}(u^2)u+2\mathcal{I}(u
v)u+\mathcal{I}(v^2)u
+\mathcal{I}(u^2)v+2\mathcal{I}(uv)v \Big]d\tau\\
\triangleq&\widetilde{\mathcal{T}}u.
\end{align*}
Using the Strichartz estimate, we have
\begin{eqnarray*}
\left\|\int_0^tK(t-\tau)\mathcal{I}(u^2)u d\tau\right\|_{L_T^\infty(
H^1)}&\lesssim &\|\mathcal{I}(u^2)u \|_{L_T^1( L^2)}.
\end{eqnarray*}
 On one hand, we make use of H\"{o}lder's inequality and Hardy-Littlewood-Sobolev inequality to deduce that
\begin{eqnarray}\label{22}
\|\mathcal{I}(u^2)u \|_{L_T^1(L^2)} &\leq&
C\|\mathcal{I}(u^2)\|_{L_T^{\frac 3
2}L^{\frac9\gamma}}\|u\|_{L_T^3L^{\frac{18}{9-2\gamma}}}\cr &\leq&
C\|u\|_{L_T^3L^{\frac{18}{9-2\gamma}}}^3 \leq CT\|u\|_{L_T^\infty
H^1}^3.
\end{eqnarray}
For the rest of terms, arguing similarly as above, it can be
obtained that
\begin{eqnarray}\label{24}
\|\mathcal{I}(uv)u \|_{L_T^1(L^2)} &\leq&
C\|u\|^2_{L_T^{\infty}L^{6}} \|v\|_{L_T^{1}L^{\frac6{7-2\gamma}}}\cr
&\leq& CT^{\frac{5-\gamma}3}\|u\|_{L_T^\infty H^1}^2\|v\|_{
X^{\alpha}},
\end{eqnarray}
\begin{eqnarray}\label{25}
\|\mathcal{I}(v^2)u \|_{L_T^1(L^2)} &\leq&
C\|u\|_{L_T^{\infty}L^{6}} \|v\|^2_{L_T^{2}L^{\frac6{4-\gamma}}}\cr
&\leq& CT^{\frac{4-\gamma}3}\|u\|_{L_T^\infty H^1}\|v\|^2_{
X^{\beta}},
\end{eqnarray}
\begin{eqnarray}\label{26}
\|\mathcal{I}(u^2)v \|_{L_T^1(L^2)}\leq
CT^{\frac{5-\gamma}3}\|u\|_{L_T^\infty H^1}^2\|v\|_{ X^{\alpha}},
\end{eqnarray}
\begin{eqnarray}\label{27}
\|\mathcal{I}(uv)v \|_{L_T^1(L^2)}\leq
CT^{\frac{4-\gamma}3}\|u\|_{L_T^\infty H^1}\|v\|^2_{ X^{\beta}}.
\end{eqnarray}
A combination of \eqref{22}, \eqref{24}-\eqref{27} and the
Strichartz estimate in Proposition \ref{pro1} lead to the estimate
\begin{eqnarray*}
\|u\|_{L_T^\infty( H^1)}&\lesssim &\|u_0\|_{H^1}+\|u_1\|_{L^2}
+T\|u\|^3_{L_T^\infty(
H^1)}\\
&&+T^{\frac{5-\gamma}3}\|u\|_{L_T^\infty H^1}^2\|v\|_{
X^{\alpha}}+T^{\frac{4-\gamma}3}\|u\|_{L_T^\infty H^1}\|v\|^2_{
X^{\beta}}.
\end{eqnarray*}
As long as choosing $T$ is small enough, $\widetilde{\mathcal{T}}$
is a contraction mapping in ball $B(0,2C\mathcal{E}_{\ell, 1})$. By
means of Picard's fixed point argument we have an unique solution
$u$ to \eqref{21} in $L^\infty([0,T];H^1)$. Therefore, Proposition
\ref{prop6} is proved  by  the standard argument. \vskip0.3cm

\section{Energy estimate for the low frequency part}
\vskip0.3cm In order to  extend  the local solution to a global
solution, we shall prove a prior estimate for the Hamiltonian of $u$
in this section. Let us recall the definition of Hamiltonian of $u$
defined by
$$H(u)(t)\stackrel{\mathrm{def}}{=}\left(\frac12\|u(t)\|^2_{H^1}
+\frac12\|u_t(t)\|^2_{L^2}+\frac14\int_{\R^3\times\R^3}|x-y|^{-\gamma}u^2(y,t)u^2(x,t)dydx\right)$$
Similarly we give another notation of the energy of $u$,  which is
denoted  by
$$E(u)(t)\stackrel{\mathrm{def}}{=}\frac12\|u(t)\|^2_{H^1}
+\frac12\|u_t(t)\|^2_{L^2}.$$ Let
$$H_T(u)\stackrel{\mathrm{def}}{=\!=}\sup_{t\leq T}H(u)(t),\quad
E_T(u)\stackrel{\mathrm{def}}{=\!=}\sup_{t\leq T}E(u)(t).$$  To
extend the local existence to global existence, we have to do a
number of nonlinear a priori estimates  provided  that $E_T(u)\leq
2CH(u)(0)$, see Proposition \ref{prop7} and Lemma \ref{lem6}. As a
direct consequence of the above assumption, we get  an important
relationship between $E(u)$ and $\mathcal{E}_s$ defined in the
introduction
\begin{align}\label{28add1} E_T(u)\lesssim
2^{2J(1-s)}(\mathcal{E}^2_s+\mathcal{E}^4_s)\lesssim 2^{2J(1-s)}.
\end{align}
In fact, it follows from Hardy-Littlewood Sobolev inequality and the
definition of $u_0$ that
$$\|(|x|^{-\gamma}\ast u_0^2)u_0^2\|_{L^1}\lesssim \|u_0\|^4_{\frac{12}{6-\gamma}}\leq
\|S_0\phi_0\|^4_{\frac{12}{6-\gamma}}+\sum_{0\leq j\leq
J}\|\triangle_j\phi_0\|^4_{\frac{12}{6-\gamma}}.$$ And then the
right hand of the above inequality can controlled that as soon as
$1>s>\frac\gamma4$ by utilizing Bernstein inequality
$$\|S_0\phi_0\|^4_{L^2}+\sum_{0\leq j\leq J}2^{j4(\frac\gamma4-s)}2^{j4s}\|\triangle_j\phi_0\|^4_{L^2}\lesssim 2^{2J(1-s)}\mathcal{E}^4_s.$$
From now on, we assume \eqref{28add1} to in our subsequence proof.
\begin{proposition}\label{prop7}\quad  Assume that $(u_{0},u_{1})\in
H^1\times L^2$, then the following estimate holds  for $s_0,\alpha,
\beta$ defined in Proposition \ref{prop3} and Proposition
\ref{prop6},
\begin{align*}
  H_T(u)\lesssim&
  H(u)(0)+T^{\frac{4-\gamma}{3}}2^{-2J(s-\beta)}E_T(u)
  +T^{\frac{5-\gamma}3}
  2^{-J(4s-\alpha-2s_0-1)}E_T(u)\cr
   &+\big(T^{\frac12+\frac1{r_1}}2^{-2J[s-(\frac\gamma2-\frac34+\frac1{2r_1})]}+
T^{\frac12+\frac1{r_2}}2^{-2J[s-(\frac\gamma2-\frac34+\frac1{2r_2})]}+T2^{-2J(s-\frac12)}\big)E_T(u)
\end{align*}
for $\max\{2,\frac1{3-\gamma}\}<r_1<\frac2{3-\gamma}$ and
$\frac4{\gamma-2}\leq r_2<\infty$.
\end{proposition}
 {\bf Proof.} Multiplying \eqref{21} by $\partial_tu$ and
integrating over $x$ and $t$, we have
\begin{align*}
&\frac12\left(\|u(t)\|^2_{H^1}
+\|u_t(t)\|^2_{L^2}\right)+\frac12\int_{\R^3\times\R^3}|x-y|^{-\gamma}u^2(y,t)u^2(x,t)dydx\\
\leq &\frac12\left(\|u_0\|^2_{H^1}
+\|u_1\|^2_{L^2}\right)+\frac12\int_{\R^3\times\R^3}|x-y|^{-\gamma}u_0^2(y)u_0^2(x)dydx\\
&+\left|\int_{0}^t\int_{\R^3}
\mathcal{I}(v^2)(x,\tau)u(x,\tau)\partial_\tau u(x,\tau)dx
d\tau\right|\\
&+2\left|\int_{0}^t\int_{\R^3}
\mathcal{I}(uv)(x,\tau)v(x,\tau)\partial_\tau u(x,\tau)dx
d\tau\right|\\
&+\left|\int_{0}^t\int_{\R^3}
\mathcal{I}(u^2)(x,\tau)v(x,\tau)\partial_\tau u(x,\tau)dx
d\tau\right|\\
&+2\left|\int_{0}^t\int_{\R^3}
\mathcal{I}(uv)(x,\tau)u(x,\tau)\partial_\tau u(x,\tau)dx
d\tau\right|.
\end{align*}
By taking the supermum over $t\leq T$, we have
\begin{align}
H_T(u)\lesssim & H(u)(0)+\left\| \mathcal{I}(v^2)u\partial_t
u\right\|_{L^1_TL^1}+\left\| \mathcal{I}(uv)v\partial_t
u\right\|_{L^1_TL^1}\nonumber\\
&+\Big|\int_0^T\int_{\R^3} \mathcal{I}(u^2)v\partial_t
u\mathrm{d}x\mathrm{d}t\Big|+\Big|\int_0^T\int_{\R^3}
\mathcal{I}(uv)u\partial_t
u\mathrm{d}x\mathrm{d}t\Big|\nonumber\\
\stackrel{\mathrm{def}}{=}&H(u)(0)+I+II+III+IV.\label{28add3}
\end{align}
The proof is broken down into the following several steps.

(i) Firstly, we estimate $I$ and $II$. Making a similarly argument
as \eqref{25} in the proof of Proposition \ref{prop6}, it can be
obtained that
\begin{align*}
I\leq\|\mathcal{I}(v^2)u\|_{L^1_TL^2}\|u_t\|_{L^\infty_TL^2}\leq
T^{\frac{4-\gamma}3}E_T(u)\|v\|^2_{ X^{\beta}},
\end{align*}
and then keeping in mind $v$ has been estimated in Proposition
\ref{prop3}, this together with \eqref{3} yields that
\begin{align}\label{29}
I\leq T^{\frac{4-\gamma}3}E_T(u)\mathcal{E}^2_{h,\beta} \leq
T^{\frac{4-\gamma}3}E_T(u)2^{-2J(s-\beta)}\mathcal{E}_s^2.
\end{align}
Arguing similarly, we easily get
\begin{align}\label{30}
II \leq T^{\frac{4-\gamma}3}E_T(u)2^{-2J(s-\beta)}\mathcal{E}_s^2.
\end{align}
(ii) Secondly, we estimate the terms $III$ and $IV$. As mentioned in
the introduction, one can get the same type of estimate as above for
the terms $I$ and $II$, but that will lead to
$s>\frac\alpha2+\frac12$, which is worse than the exponent given in
the Theorem \ref{the1}. To improve the lower bound on $s$, we have
to utilize more precise estimate on $III$ and $IV$.

We first split $III$ and $IV$ into two different pieces,
respectively.  One can write
$$v=v_F+B(v,v,v),$$
where $v_F$ is its free part and the other one comes from nonlinear
term. For the nonlinear part, it follows from \eqref{7a} that
\begin{align*}
\|B(v,v,v)\|_{X^{\alpha}} \leq \|v\|_{X^{\alpha}}\|v\|^2_{X^{s_0}}.
\end{align*}
This along with \eqref{26}, one can see that
\begin{equation*}
\begin{split}
\|\mathcal{I}(u^2)B(v,v,v)u_t\|_{L^1_TL^1} &\leq
\|\mathcal{I}(u^2)B(v,v,v)\|_{L^1_TL^2}\|u_t\|_{L^\infty_TL^2}\\
&\leq T^{\frac{5-\gamma}3}\|u\|^2_{L^\infty_TH^1}
\|B(v,v,v)\|_{X^{\alpha}}\|u_t\|_{L^\infty_TL^2}\\
&\leq T^{\frac{5-\gamma}3}E_T(u)^{\frac32}
\|v\|_{X^{\alpha}}\|v\|^2_{X^{s_0}}\\
\end{split}
\end{equation*}
Moreover,  we get by \eqref{3},
\begin{align}\label{31}
\|\mathcal{I}(u^2)B(v,v,v)u_t\|_{L^1_TL^1} &\leq
T^{\frac{5-\gamma}3}E^{\frac32}_T(u)
2^{-J(3s-\alpha-2s_0)}\mathcal{E}_s^3.
\end{align}
By the same way  as leading to \eqref{31}, we easily infer that
\begin{align}\label{32}
\|\mathcal{I}(uB(v,v,v))uu_t\|_{L^1_TL^1} &\leq
T^{\frac{5-\gamma}3}E^{\frac32}_T(u)
2^{-J(3s-\alpha-2s_0)}\mathcal{E}_s^3.
\end{align}
Thus, it is sufficient to estimate these terms including free part
$v_F$ since \eqref{31} and \eqref{32}. The following lemma gives
estimates for the nonlinearity including free part $v_F$.
\begin{lemma}\label{lem6}
Let $v_F$ be a solution of the free Klein-Gordon equation, and $u$
be such that $E_T(u)\lesssim 2^{2J(1-s)}$. Then, for
$\max\{2,\frac1{3-\gamma}\}<r_1<\frac2{3-\gamma}$ and
$\frac4{\gamma-2}\leq r_2<\infty$
\begin{align}
\label{33}\Big|\int_0^T\int_{\R^3}\mathcal{I}(u^2)v_Fu_t\mathrm{d}x\mathrm{d}t\Big|
&\lesssim
\Big(T^{\frac12+\frac1{r_1}}2^{-2J[s-(\frac\gamma2-\frac34+\frac1{2r_1})]}\cr&+
T^{\frac12+\frac1{r_2}}2^{-2J[s-(\frac\gamma2-\frac34+\frac1{2r_2})]}+T2^{-2J(s-\frac12)}\Big)
E_T{(u)},\\
\label{34}\Big|\int_0^T\int_{\R^3}\mathcal{I}(uv_F)uu_t\mathrm{d}x\mathrm{d}t\Big|
&\lesssim
T^{\frac12+\frac1{r_2}}2^{-2J[s-(\frac\gamma2-\frac34+\frac1{2r_2})]}
E_T{(u)}.
\end{align}
\end{lemma}
Hence these together with \eqref{31}-\eqref{32} yield that
\begin{align}
\label{35}III+IV \lesssim T^{\frac{5-\gamma}3}E_T(u)
2^{-J[4s-\alpha-2s_0-1]}\mathcal{E}_s^4+\Big(T^{\frac12+\frac1{r_1}}2^{-2J[s-(\frac\gamma2-\frac34+\frac1{2r_1})]}\cr+
T^{\frac12+\frac1{r_2}}2^{-2J[s-(\frac\gamma2-\frac34+\frac1{2r_2})]}+T2^{-2J(s-\frac12)}\Big)E_T{(u)}.
\end{align}
Therefore, we complete the proof of Proposition \ref{prop7} provided
that we had proved Lemma \ref{lem6},  whose proof is postponed in
the last section.

\vskip0.3cm
\section{Proof of Theorem \ref{the1}}

\vskip0.3cm Since the Cauchy problem \eqref{1} is split into
equation \eqref{14} which is globally well-posed by choosing $J$
enough to make $\mathcal{E}_{h,s_0}<\varepsilon_0$ and equation
\eqref{21} which is locally well-posed (see Proposition \ref{prop3}
and Proposition \ref{prop6}),  we have to show that the local
solution to equation \eqref{21} can be extended globally.

Let us denote $T^*_J$ the maximum time of existence in Proposition
\ref{prop6}. Theorem \ref{the1} will be proved if
$$\lim_{J\rightarrow+\infty}T^*_J=+\infty.$$
Let us consider $T_J$ the supremum of the $T<T_J^*$ such that
$$E_T(u)\leq2CH(u)(0).$$ Thus, for any $T<T_J$, Proposition
\ref{prop7} gives us that
\begin{align*}
  E_T(u)\leq
  H(u)(0)\bigg(&C+C_1T^{\frac{4-\gamma}{3}}2^{-2J(s-\beta)}\mathcal{E}_s^2\\
  &+C_2T^{\frac{5-\gamma}3} 2^{-J(4s-\alpha-2s_0-1)}\mathcal{E}_s^4
   +C_3T^{\frac{1}2+\frac1{r_1}}
   2^{-2J[s-(\frac\gamma2-\frac34+\frac1{2r_1})]}\mathcal{E}_s^2\\&+C_4T
   2^{-2J(s-\frac12)}\mathcal{E}_s^2+C_5T^{\frac12+\frac1{r_2}}2^{-2J[s-(\frac\gamma2-\frac34+\frac1{2r_2})]}
\mathcal{E}^2_{s}\bigg).
\end{align*}
By the assumption of Theorem \ref{the1} $s>\frac\gamma4$, one easily
verifies that
\begin{align*}
s>\max\left\{\beta, \frac\alpha4+\frac{s_0}2+\frac14, \frac12,
\frac\gamma2-\frac34+\frac1{2r_1},\frac\gamma2-\frac34+\frac1{2r_2}
\right\}
\end{align*}
if choosing $r_1$ sufficiently close to $\frac2{3-\gamma}$ and $r_2$
large enough. We infer that $T_J\geq \widetilde{T}_J$ if we choose
$\widetilde{T}_J$ such that
\begin{align*}
\widetilde{T}_J\stackrel{\mathrm{def}}{=}\min\Bigg\{&
\left(\frac{2^{2J(s-\beta)}}{5C_1\mathcal{E}^2_s}\right)^{\frac3{4-\gamma}},
\left(\frac{2^{4J(s-\frac1{4}\alpha-\frac{s_0}2-\frac14)}}{5C_2\mathcal{E}^4_s}\right)
^{\frac3{5-\gamma}},
\left(\frac{2^{2J[s-(\frac\gamma2-\frac34+\frac1{2r_1})]}}{5C_3\mathcal{E}^{2}_s}\right)^{\frac{2r_1}{r_1+2}},\cr&
\frac{2^{2J(s-\frac12)}}{5C_4\mathcal{E}^{2}_s},\quad
\left(\frac{2^{2J[s-(\frac\gamma2-\frac34+\frac1{2r_2})]}}{5C_5\mathcal{E}^{2}_s}\right)^{\frac{2r_2}{r_2+2}}\Bigg\}.
\end{align*}
By the definition of ${T}_J$, we get $T^*_J\geq \widetilde{T}_J$.
Obviously, $\widetilde{T}_J$ tends to infinity when $J$ tend to
infinity. This completes the proof of  Theorem \ref{the1}.

\section{Proof of Lemma \ref{lem6}}
In order to make conveniently use of the precise Strichartz estimate
on which mostly the following proof relies, we begin this section by
introducing a family of balls of center
$(\xi^{j,k}_\nu)_{\nu\in\Lambda_{j,k}}$ of radius $2^k$ and a
function $\chi\in\mathcal{C}_c^\infty(B(0,1))$ such that for $j\geq
0$
$$\forall \xi\in2^j\mathcal{C},\sum_{\nu\in\Lambda_{j,k}}\chi(2^{-k}(\xi-\xi_\nu^{j,k}))=1
 \quad and \quad C_0^{-1}\leq\sum_{\nu\in\Lambda_{j,k}}\chi^2(2^{-k}(\xi-\xi_\nu^{j,k}))\leq C_0.$$
Let us define that, for some constant $c$
\begin{eqnarray*}
\triangle_{j,k}^\nu a
&\stackrel{\mathrm{def}}{=\!=}&\mathcal{F}^{-1}\big(\big(\varphi(2^{-j}\xi)\chi(2^{-k}(\xi-\xi_\nu^{j,k}))\big)\hat
a(\xi)\big),\\
\widetilde{\triangle_{j,k}^\nu} a
&\stackrel{\mathrm{def}}{=\!=}&\mathcal{F}^{-1}\big(\big(\widetilde{\varphi}(2^{-j}\xi)\chi(c2^{-k}(\xi+\xi_\nu^{j,k}))\big)\hat
a(\xi)\big).
\end{eqnarray*}
As the support of the Fourier transform of a product belongs to  the
sum of the support of each Fourier transform, we have
\begin{eqnarray*}
\triangle_j a=\sum_{\nu\in\Lambda_{j,k}}\triangle_{j,k}^\nu a, \quad
\triangle_j b=\sum_{\nu'\in\Lambda_{j,k}}\triangle_{j,k}^{\nu'} b.
\end{eqnarray*}
In view of this fact that if $k\leq j-2$
$$\triangle_k\sum_{\nu,\nu'\in\Lambda_{j,k}}\triangle_{j,k}^\nu
a\triangle_{j,k}^{\nu'} b$$ is vanish when $\xi_\nu^{j,k}$ is close
to $\xi_{\nu'}^{j,k}$, without loss of generality, we can write
\begin{eqnarray}\label{add1}\triangle_k(
\triangle_ja\triangle_jb)\approx\triangle_k\sum_{\nu\in\Lambda_{j,k}}\triangle_{j,k}^\nu
a\widetilde{\triangle_{j,k}^{\nu}} b.\end{eqnarray} For the sake of
convenience, we also fix the notation in this section that, for
$0\neq f(t,x)\in L^2_TL^2$
\begin{align*} c_k=2^{k\sigma} \bigg(\|\triangle_k
v_0\|_{L^2}+2^{-k}\|\triangle_k
v_1\|_{L^2}\bigg){\mathcal{E}^{-1}_{h,\sigma}},\quad
\tilde{c}_k=\frac{\|\triangle_kf\|_{L^2_TL^{2}_x}}{\|f\|_{L^2_TL^2}}
\end{align*}
with $\sigma=1/2+1/r$ for $2\leq r<\infty$.

{\bf Proof of Lemma \ref{lem6}}. We first prove \eqref{33}. In view
of the fact that $\widehat{v_F}$ only has high frequencies, Bony's
decomposition implies that there exists constant $N_0$ such that
\begin{eqnarray}\label{7.0}
\mathcal{I}(u^2)v_Fu_t =\displaystyle\sum_{j\geq J-N_0}
S_{j+2}{v_F}\triangle_j\mathcal{I}(u^2)u_t+ \displaystyle\sum_{j\geq
J-N_0} S_{j-1}{\mathcal{I}(u^2)}\triangle_jv_Fu_t.
\end{eqnarray}
Since the negative derivative $\mathcal{I}$ acts on the high
frequency for the former term while on the low frequency for the
latter one, the first term is much better than the second one. We
shall estimate the first term by using merely the H\"{o}lder
inequality, Bernstein inequality and classical Strichartz estimates.
Firstly, we see that, for $2\leq r<\infty$
\begin{align*}
\sum_{j\geq J-N_0}\|
S_{j+2}{v_F}\triangle_j\mathcal{I}(u^2)\|_{L_x^2}&\lesssim\sum_{j\geq
J-N_0}\sum_{j^\prime\leq
j+1}\|\triangle_{j^\prime}{v_F}\|_{L_x^{\infty}}
\|\triangle_j\mathcal{I}(u^2)\|_{L_x^{2}}\cr &\lesssim \sum_{j\geq
J-N_0}\sum_{j^\prime\leq j+1}2^{j^\prime\frac3
r}\|\triangle_{j^\prime}{v_F}\|_{L^r}2^{j(\gamma-\frac72)}\|u\|^2_{L_T^\infty
H^1}.
\end{align*}
 Bernstein inequality and \eqref{11} in Proposition
\ref{pro1} with $\frac1p+\frac1r=\frac12$ for all $2\leq r<\infty$
imply that
\begin{align*}
\big\|\sum_{j\geq J-N_0}&\|
S_{j+2}{v_F}\triangle_j\mathcal{I}(u^2)\|_{L_x^2}\big\|_{L^1_T}
\lesssim T^{1-\frac1p}\|u\|^2_{L_T^\infty H^1}\sum_{j\geq
J-N_0}\sum_{j^\prime\leq j+1}2^{j^\prime\frac3
r}\|\triangle_{j^\prime}{v_F}\|_{L^p_TL^r}2^{j(\gamma-\frac72)}\cr
&\lesssim T^{1-\frac1p}\|u\|^2_{L_T^\infty H^1}\sum_{j\geq
J-N_0}2^{j(\gamma-\frac 7 2)}\sum_{j^\prime\leq
j+1}2^{j^\prime(\frac32-\frac1p-\sigma)}2^{j^\prime\sigma}
\left(\|\triangle_{j^\prime}{v_0}\|_{L^2}+2^{-j^\prime}\|\triangle_{j^\prime}{v_1}\|_{L^2}\right)
.
\end{align*}
The right hand of the above inequality can be controlled by
\begin{equation*}
\begin{split}
T^{1-\frac1p}\|u\|^2_{L_T^\infty H^1}\sum_{j\geq
J-N_0}2^{j(\gamma-\frac72)}\sum_{j^\prime\leq
j+1}2^{\frac{j'}2}c_{j'}\mathcal{E}_{h,\sigma}
\end{split}
\end{equation*}
and moreover it follows from \eqref{3}, the definition of
$\mathcal{E}_{h,\sigma}$ and Sobolev embedding that
\begin{equation}\label{37}
\begin{split}
\big\|\sum_{j\geq J-N_0}\|
S_{j+2}{v_F}\triangle_j\mathcal{I}(u^2)u_t\|_{L_x^1}\big\|_{L^1_T}&\lesssim
T^{1-\frac1p}\mathcal{E}_{h,\sigma}\|u\|^2_{L_T^\infty H^1}
\sum_{j\geq J-N_0}2^{j(\gamma-3)}\|u_t\|_{L_T^\infty L^2}\cr
&\lesssim
T^{\frac12+\frac1r}2^{-2J[s-(\frac\gamma2-\frac34+\frac1{2r})]}\mathcal{E}^2_{s}\|u\|^2_{L^\infty_T(H^1)},
\end{split}
\end{equation}
for $\frac{4}{\gamma-2}\leq r<\infty$.

\vskip 0.1in Let us estimate the second term in \eqref{7.0} by the
precise Strichartz estimates. Since this term contains that the
negative derivative acts on the low frequency part $S_{j-1}(u^2)$,
it leads to our new parameter $r<\frac2{3-\gamma}$ by some technique
difficulties. Noting that Fourier-Plancherel formula and H\"older's
inequality, we can see that
\begin{equation}\label{7.4}
\begin{split}
\sum_{j\geq J-N_0} &\int_0^T\int_{\R^3}S_{j-1}{\mathcal{I}(u^2)}
\triangle_jv_Fu_t\mathrm{d}x\mathrm{d}t\lesssim\int \sum_{j\geq
J-N_0}\sum_{-1\leq k\leq j-2}\triangle_k{\mathcal{I}(u^2)}
\triangle_jv_F \triangle_ju_tdxdt\cr &\approx\sum_{k\geq-1}\int
\triangle_k{\mathcal{I}(u^2)}\triangle_k\sum_{k\leq j-2, J-N_0\leq
j}\big( \triangle_jv_F
\triangle_ju_t\big)dxdt\cr&\lesssim\|u^2\|_{L^\infty_T(
B^{\frac12}_{2,1})}\int_0^T
\sup_{k\geq-1}2^{k(\gamma-\frac72)}\big\|\triangle_k\sum_{k\leq j-2,
J-N_0\leq j}\big( \triangle_jv_F \triangle_ju_t\big)\big\|_{L^2}dt.
\end{split}
\end{equation}
On one hand, we have
\begin{align*}
\int_0^T \|\triangle_{-1}\sum_{J-N_0\leq j}\big( \triangle_jv_F
\triangle_ju_t\big)\|_{L^2}dt&\lesssim \sum_{j\geq
J-N_0}\|\triangle_jv_F \triangle_ju_t\|_{L^1_TL^{1}}\cr &\leq
T^{\frac12}\sum_{j\geq J-N_0}\|\triangle_jv_F\|_{L^\infty_TL^{2}}
\|\triangle_ju_t\|_{L^2_TL^{2}}\cr &\leq T^{\frac12}\sum_{j\geq
J-N_0}2^{-sj}c_j\tilde{c}_j\mathcal{E}_{h,s}\|u_t\|_{L^2_{t,x}}.
\end{align*}
If \eqref{7.4} is controlled by the term  at $k=-1$, we can see that
\begin{align}\label{38} \sum_{j\geq J-N_0}\|
S_{j-1}{\mathcal{I}(u^2)} \triangle_jv_F\|_{L^2_TL^2}&\lesssim
T^{\frac12}2^{-2J(s-\frac12)}\mathcal{E}^2_{s}\|u\|_{L^\infty_T(H^1)}.
\end{align}
On the other hand, one denotes $ g_k:=\triangle_k\sum\limits_{k\leq
j-2}\big( \triangle_jv_F \triangle_ju_t\big)$ to estimate
$$\sum_{k\geq 0}2^{k(\gamma-\frac72+\frac3r)}\|g_k\|_{L^1_TL^{\frac{2r}{r+2}}}.$$
Let us write that
\begin{align*}
g_k&= \sum_{k\leq
j-2}\triangle_k\sum_{\nu\in\Lambda_{j,k}}\triangle_{j,k}^\nu
v_F\triangle_ju_t.
\end{align*}
As the support of the Fourier transform of a product is included in
the sum of the support of each Fourier transform, we obtain
\begin{align*}
g_k&= \sum_{k\leq
j-2}\triangle_k\sum_{\nu\in\Lambda_{j,k}}\triangle_{j,k}^\nu
v_F\widetilde{\triangle_{j,k}^\nu}u_t,
\end{align*}
as well as in \eqref{add1}. Using H\"older inequality, we get
\begin{align*}
\|g_k\|_{L^{\frac{2r}{r+2}}}&\leq \sum_{k\leq
j-2}\sum_{\nu\in\Lambda_{j,k}}\|\triangle_{j,k}^\nu
v_F\|_{L^r}\|\widetilde{\triangle_{j,k}^\nu}u_t\|_{L^2}
\end{align*}
and the Cauchy-Schwarz inequality and the $L^2$ quasi-orthogonality
properties yield that
\begin{align}\label{39}
\|g_k\|_{L^{\frac{2r}{r+2}}}&\leq \sum_{k\leq
j-2}\big(\sum_{\nu\in\Lambda_{j,k}}\|\triangle_{j,k}^\nu
v_F\|^2_{L^r}\big)^{\frac12}\big(\sum_{\nu\in\Lambda_{j,k}}\|\widetilde{\triangle_{j,k}^\nu}
u_t\|^2_{L^2}\big)^{\frac12}\cr &\leq \sum_{k\leq
j-2}\big(\sum_{\nu\in\Lambda_{j,k}}\|\triangle_{j,k}^\nu
v_F\|^2_{L^r}\big)^{\frac12}\|\triangle_{j} u_t\|_{L^2}.
\end{align}
Precise Strichartz estimate implies that, for $\frac1 p+\frac1
r=\frac12$ with $2\leq r<\infty$,
\begin{align*}\|g_k\|_{L^1_T(L^{\frac{2r}{r+2}})}&\leq
T^{\frac12-\frac1p} \sum_{0\leq k\leq
j-2}2^{(k-j)(\frac12-\frac1r)}2^{j(\frac32-\frac3r-\frac1p)}\cr&\,\,\,\,\,\,\,\,\,
\times \bigg(\big(\sum_{\nu\in\Lambda_{j,k}}\|\triangle_{j,k}^\nu
v_0\|^2_{L^2}\big)^{\frac12}+
2^{-j}\big(\sum_{\nu\in\Lambda_{j,k}}\|\triangle_{j,k}^\nu
v_1\|^2_{L^2}\big)^{\frac12}\bigg)\|\triangle_ju_t\|_{L^2_{t,x}}
\end{align*}
and observe the quasi-orthogonality properties again, this can be
dominated by
\begin{align*} T^{\frac12-\frac1p} \sum_{0\leq k\leq
j-2}2^{(k-j)(\frac12-\frac1r)}2^{j(\frac32-\frac3r-\frac1p)}
\bigg(\|\triangle_j v_0\|_{L^2}+ 2^{-j}\|\triangle_j
v_1\|_{L^2}\bigg)\|\triangle_ju_t\|_{L^2_{t,x}}.
\end{align*}
Keeping the definitions of $\mathcal{E}_{h,\sigma}$ and $c_j$ in
mind, one can see that
\begin{align*} \|g_k\|_{L^1_T(L^{\frac{2r}{r+2}})}\lesssim
T^{\frac12-\frac1p}2^{k(\frac12-\frac1r)} \sum_{k\leq j-2}2^{-\frac
{2j}r}c_j\tilde{c}_j \mathcal{E}_{h,\sigma}E^{\frac12}_T(u)\lesssim
T^{\frac12-\frac1p}2^{k(\frac12-\frac3r)}
\mathcal{E}_{h,\sigma}E^{\frac12}_T(u).
\end{align*}
Therefore, we get that
\begin{align*} \sum_{k\geq0}2^{k(\gamma-\frac72+\frac3r)}\|g_k\|_{L^1_T(L^{\frac{2r}{r+2}})}
&\lesssim T^{\frac12-\frac1p}\sum_{k\geq0}2^{k(\gamma-3)}
\mathcal{E}_{h,\sigma}E^{\frac12}_T(u)
\end{align*}
which implies nothing but
\begin{align}\label{7.7}
&\Big|\int_0^T\int_{\R^3}\sum_{j\geq J-N_0}S_{j-1}{\mathcal{I}(u^2)}
\triangle_jv_Fu_t\mathrm{d}x\mathrm{d}t\Big|\lesssim
T^{1-\frac1p}\mathcal{E}_{h,\sigma}E^{\frac32}_T(u).
\end{align}
Finally, we get that, for $\frac{4}{\gamma-2}\leq r<\infty$
\begin{align*}
&\Big|\int_0^T\int_{\R^3}\sum_{j\geq J-N_0}S_{j-1}{\mathcal{I}(u^2)}
\triangle_jv_Fu_t\mathrm{d}x\mathrm{d}t\Big|\lesssim
T^{\frac12+\frac1r}2^{-2J[s-(\frac34+\frac1{2r})]}\mathcal{E}^2_{s}E_T(u).
\end{align*}
However, although the $r$ ranges $\frac{4}{\gamma-2}\leq r<\infty$,
the above estimate still needs $s>\frac34$ to continue our proof. If
we only consider the high frequency $k\geq J$, the \eqref{7.7} can
be modified by
\begin{align}\label{7.8}
&\Big|\int_0^T\int_{\R^3}\sum_{j\geq J-N_0}S_{j-1}{\mathcal{I}(u^2)}
\triangle_jv_Fu_t\mathrm{d}x\mathrm{d}t\Big|\lesssim
T^{1-\frac1p}2^{J(\gamma-3)}\mathcal{E}_{h,\sigma}E^{\frac32}_T(u)
\end{align}
and then we can obtain a better result
\begin{align*}
&\Big|\int_0^T\int_{\R^3}\sum_{j\geq J-N_0}S_{j-1}{\mathcal{I}(u^2)}
\triangle_jv_Fu_t\mathrm{d}x\mathrm{d}t\Big|\lesssim
T^{\frac12+\frac1r}2^{-2J[s-(\frac\gamma2-\frac34+\frac1{2r})]}\mathcal{E}^2_{s}E_T(u),
\end{align*}
which implies the bad influence comes from the low frequency part
and this is consist of the effect of negative derivative acts on the
low frequency. But if we choose
$\tilde{\sigma}=\gamma-\frac52+\frac1r$ instead of $\sigma$, we can
improve \eqref{7.8}, at cost of restricting $r$ such that
$\max\{2,\frac{1}{3-\gamma}\}<r<\frac{2}{3-\gamma}$ while not $2\leq
r<\infty$. Now we turn to details. It follows from similar argument
that
\begin{align*} \|g_k\|_{L^1_T(L^{\frac{2r}{r+2}})}&\lesssim
T^{\frac12-\frac1p}2^{k(\frac12-\frac1r)} \sum_{k\leq
j-2}2^{-j(\gamma-3+\frac 2r)}c_j\tilde{c}_j
\mathcal{E}_{h,\tilde{\sigma}}E^{\frac12}_T(u)
\end{align*}
where $\tilde{\sigma}=\gamma-\frac52+\frac1r$ with
$\frac1r<3-\gamma<\frac2r$. We get
\begin{align*} \sum_{k\geq0}2^{k(\gamma-\frac72+\frac3r)}\|g_k\|_{L^1_T(L^{\frac{2r}{r+2}})}
&\lesssim T^{\frac12-\frac1p}\sum_{k\geq0}\sum_{k\leq
j-2}2^{(k-j)(\gamma-3+\frac 2r)}c_j\tilde{c}_j
\mathcal{E}_{h,\tilde{\sigma}}E^{\frac12}_T(u)
\end{align*}
which implies nothing but
\begin{align*}
&\Big|\int_0^T\int_{\R^3}\sum_{j\geq J-N_0}S_{j-1}{\mathcal{I}(u^2)}
\triangle_jv_Fu_t\mathrm{d}x\mathrm{d}t\Big|\lesssim
T^{1-\frac1p}\mathcal{E}_{h,\tilde{\sigma}}E^{\frac32}_T(u)
\end{align*}
by Young's inequality. Note that $\tilde{\sigma}\leq \frac\gamma4<s$
when $r$ sufficiently closes to $\frac2{3-\gamma}$, therefore
\begin{align*}
&\Big|\int_0^T\int_{\R^3}\sum_{j\geq J-N_0}S_{j-1}{\mathcal{I}(u^2)}
\triangle_jv_Fu_t\mathrm{d}x\mathrm{d}t\Big|\lesssim
T^{\frac12+\frac1r}2^{-2J(s-\frac{\tilde{\sigma}+1}2)}\mathcal{E}^2_{s}E_T(u).
\end{align*}
Combining this with \eqref{37} and \eqref{38}, we  complete the
proof of \eqref{33} by obtaining
\begin{align*}
&\Big|\int_0^T\int_{\R^3}\sum_{j\geq J-N_0}S_{j-1}{\mathcal{I}(u^2)}
\triangle_jv_Fu_t\mathrm{d}x\mathrm{d}t\Big|\cr&\lesssim
\Big(T^{\frac12+\frac1{r_1}}2^{-2J[s-(\frac\gamma2-\frac34+\frac1{2r_1})]}+
T^{\frac12+\frac1{r_2}}2^{-2J[s-(\frac\gamma2-\frac34+\frac1{2r_2})]}+T2^{-2J(s-\frac12)}\Big)
\mathcal{E}^2_{s}E_T{(u)}
\end{align*}
with $\max\{2,\frac{1}{3-\gamma}\}<r_1<\frac{2}{3-\gamma}$ and
$\frac{4}{\gamma-2}\leq r_2<\infty$.

\vskip 0.1in We secondly prove \eqref{34} which is different from
\eqref{33}. To this end, we need to make Bony's decomposition more
than once and establish a commutator estimate, which helps us to
complete our proof. In view of the fact that $\widehat{v_F}$ only
has high frequencies again, it follows from Bony's decomposition
that there exists $N_0$ such that
\begin{eqnarray}\label{}
\mathcal{I}(uv_F)uu_t =\displaystyle\sum_{j\geq
J-N_0}\mathcal{I}\big( S_{j+2}{v_F}\triangle_ju\big)uu_t+
\displaystyle\sum_{j\geq J-N_0}\mathcal{I}\big(\triangle_jv_F
S_{j-1}{u}\big)uu_t\stackrel{\mathrm{def}}{=\!=}I+II.
\end{eqnarray}
In order to estimate the term $I$, we split it into two pieces with
$N_1\gg N_0>0$
\begin{align*}
I&=\sum_{j\geq J-N_0}\sum_{k} uu_t\triangle_k\mathcal{I}\big(
S_{j+2}{v_F}\triangle_ju\big)\cr&=\sum_{j\geq J-N_0}\sum_{k\leq
J-N_1}uu_t\triangle_k\mathcal{I}\big(
S_{j+2}{v_F}\triangle_ju\big)+\sum_{j\geq J-N_0}\sum_{k\geq
J-N_1}uu_t\triangle_k\mathcal{I}\big(
S_{j+2}{v_F}\triangle_ju\big)\cr&\stackrel{\mathrm{def}}{=\!=}I_1+I_2.
\end{align*}
The estimate of $I_1$ is broken down into the following two cases.
\vskip 0.1in {\bf Case 1.} $2<\gamma\leq\frac52$ \vskip 0.1in In
this case, to our purpose, we obtain the following coarse estimate
by H\"{o}lder's inequality
\begin{align*}
\|I_1\|_{L^1_x}&\lesssim\sum_{j\geq J-N_0}\sum_{k\leq
J-N_1}\|\triangle_k\mathcal{I}\big(
S_{j+2}{v_F}\triangle_ju\big)\|_{L^3}\|u\|_{L^\infty_T
L^6}\|u_t\|_{L^\infty_T L^2}\cr&\lesssim\sum_{j\geq
J-N_0}\sum_{k\leq J-N_1}2^{k(\gamma-2)}\|\triangle_k\big(
S_{j+2}{v_F}\triangle_ju\big)\|_{L^{\frac32}}E_T(u)\cr&\lesssim\sum_{j\geq
J-N_0}\sum_{k\leq J-N_1}2^{k(\gamma-2)}\|
S_{j+2}{v_F}\|_{L^6}\|\triangle_ju\|_{L^2}E_T(u)\cr&\lesssim
\sum_{k\leq J-N_1}2^{k(\gamma-2)}\sum_{j\geq
J-N_0}2^{-j}2^j\|\triangle_ju\|_{L^2}\sum_{j^\prime\leq j}\|
\triangle_{j^\prime}{v_F}\|_{ L^6}E_T(u).
\end{align*}

Choosing $(p,r)$ such that $\frac1p+\frac1r=\frac12$ with $2\leq
r\leq6$, the Strichartz estimate yeilds
\begin{align*}
\|I_1\|_{L^1_TL^1_x}\lesssim T^{1-\frac1p}\sum_{k\leq
J-N_1}2^{k(\gamma-2)}\sum_{j\geq J-N_0}&2^{-j}\sum_{ j^\prime\leq
j}2^{j^\prime(\frac3r-\frac36)}2^{j^\prime(\frac32-\frac3r-\frac1p)}
\cr&\times\big(\| \triangle_{j'}{v_0}\|_{L^2}+2^{-j^\prime}\|
\triangle_{j'}{v_1}\|_{L^2}\big)E^{\frac32}_T(u).
\end{align*}
Arguing similarly as before it yields that
\begin{align*}
\|I_1\|_{L^1_TL^1_x}&\lesssim T^{\frac12+\frac1r}\sum_{k\leq
J-N_1}2^{k(\gamma-2)}\sum_{j\geq J-N_0}2^{-j}\sum_{ j^\prime\leq
j}2^{\frac{j^\prime}r}c_{j'}\mathcal{E}_{h,1/2}
E^{\frac32}_T(u)\cr&\lesssim
T^{\frac12+\frac1r}2^{-2J[s-(\frac\gamma2-\frac34+\frac1{2r})]}
\mathcal{E}^2_{s}E_T(u)
\end{align*}
with $2\leq r\leq6$. If choose $r=6$, one can easily check that
$\frac\gamma 4>\frac\gamma2-\frac34+\frac1{2r}$ when
$2<\gamma\leq\frac52$. Although this result is enough for us to
prove the main theorem, we want to improve the result for this term
by loosen the upper bound of $r$ from $6$ to $ \infty$ through the
precise Strichartz estimate. Arguing similarly as before, we have
\begin{align*}
\|I_1\|_{L^1_x}&\lesssim\sum_{j\geq J-N_0}\sum_{k\leq
J-N_1}\|\triangle_k\mathcal{I}\big(
S_{j+2}{v_F}\triangle_ju\big)\|_{L^3}\|u\|_{L^\infty_T
L^6}\|u_t\|_{L^\infty_T L^2}\cr&\lesssim\sum_{j\geq
J-N_0}\sum_{k\leq
J-N_1}2^{k(\gamma-3)}2^{k\frac{r+6}{2r}}\|\triangle_k\big(
S_{j+2}{v_F}\triangle_ju\big)\|_{L^{\frac{2r}{r+2}}}E_T(u)
\end{align*}
Since the Fourier transform of $S_{j-1}{v_F}\triangle_ju$ was
supported in $2^j\mathcal{C}$ and $k\ll j$,
$\triangle_k(S_{j-1}{v_F}\triangle_ju)$ vanishes which implies
$\triangle_k\big( S_{j+2}{v_F}\triangle_ju\big)=\triangle_k\big(
\widetilde{\triangle_{j}}{v_F}\triangle_ju\big)$. As the support of
the Fourier transform of a product is included in the sum of the
support of each Fourier transform, we also have
\begin{align*}
\triangle_k\big(\widetilde{\triangle_{j}}{v_F}\triangle_ju\big)&=\triangle_k\big(\sum_{\nu,\nu'\in\Lambda_{j,k}}
\triangle_{j,k}^\nu{v_F}\triangle_{j,k}^{\nu'}u\big)=\triangle_k\big(\sum_{\nu\in\Lambda_{j,k}}
\triangle_{j,k}^\nu{v_F}\widetilde{\triangle_{j,k}^{\nu}}u\big)
\end{align*}
Choosing $(p,r)$ such that $\frac1p+\frac1r=\frac12$ for $2\leq
r<\infty$, it follows from the H\"{o}lder inequality and $L^2$
quasi-orthogonality properties that
\begin{align*}
\|\triangle_k\big(
S_{j+2}{v_F}\triangle_ju\big)\|_{L^1_T(L^{\frac{2r}{r+2}})}&\lesssim
\big\|\sum_{\nu\in\Lambda_{j,k}}
\|\triangle_{j,k}^\nu{v_F}\|_{L^r}\|\widetilde{\triangle_{j,k}^{\nu}}u\|_{L^2}\big\|_{L^1_T}\cr&\lesssim
T^{\frac12-\frac1p}\big(\sum_{\nu\in\Lambda_{j,k}}
\|\triangle_{j,k}^\nu{v_F}\|^2_{L^pL^r}\big)^{\frac12}\big\|\big(\sum_{\nu\in\Lambda_{j,k}}
\|\triangle_{j,k}^{\nu}u\|^2_{L^2}\big)^{\frac12}\big\|_{L^2_T}\cr&\lesssim
T^{\frac12-\frac1p} \big(\sum_{\nu\in\Lambda_{j,k}}
\|\triangle_{j,k}^\nu{v_F}\|^2_{L^pL^r}\big)^{\frac12}\big\|\triangle_{j}u\big\|_{L^2_TL^2}
\end{align*}
Then the precise Strichartz estimate yields that
\begin{align*}
\|I_1\|_{L^1_TL^1_x}&\lesssim T^{\frac12-\frac1p}\sum_{k\leq
J-N_1}2^{k(\gamma-3)}2^{k\frac{r+6}{2r}}\sum_{j\geq J-N_0
}2^{-j}2^{j}\big\|\triangle_{j}u\big\|_{L^2_TL^2}2^{(k-j)(\frac12-\frac1r)}2^{j(\frac32-\frac3r-\frac1p)}
\cr&\times\bigg(\big(\sum_{\nu\in\Lambda_{j,k}}\|
\triangle_{j,k}^\nu{v_0}\|^2_{L^2}\big)^{\frac12}+2^{-j}\big(\sum_{\nu\in\Lambda_{j,k}}\|
\triangle_{j,k}^\nu{v_1}\|^2_{L^2}\big)^{\frac12}\bigg)E_T(u).
\end{align*}
By the $ L^2$quasi-orthogonality properties, it gives that
\begin{align*}
\|I_1\|_{L^1_TL^1_x}&\lesssim T^{\frac12-\frac1p}\sum_{k\leq
J-N_1}2^{k(\gamma-3)}2^{k\frac{2r+4}{2r}}\sum_{j\geq
J-N_0}2^{-j}2^{j}\big\|\triangle_{j}u\big\|_{L^2_TL^2}2^{j(\frac12-\frac1r)}
\cr&\times\bigg(\| \triangle_{j}{v_0}\|_{L^2}+2^{-j}\|
\triangle_{j}{v_1}\|_{L^2}\bigg)E_T(u).
\end{align*}
Utilizing the technique as before yields that
\begin{align*}
\|I_1\|_{L^1_TL^1_x}&\lesssim T^{1-\frac1p}\sum_{k\leq
J-N_1}2^{k(\gamma-2+\frac2r)}\sum_{j\geq J-N_0}2^{-j(1+\frac2r)}
c_j\mathcal{E}_{h,\sigma}E^{\frac32}_T(u)\cr&\lesssim
T^{\frac12+\frac1r}2^{J(\gamma-3)}
\mathcal{E}_{h,\sigma}E^{\frac32}_T(u)\lesssim
T^{\frac12+\frac1r}2^{-2J[s-(\frac\gamma2-\frac34+\frac1{2r})]}
\mathcal{E}^2_{s}E_T(u),
\end{align*}
with $\frac4{\gamma-2}\leq r<\infty$.
\vskip 0.1in {\bf Case 2.}
$\frac52<\gamma<3$ \vskip 0.1in
In the this case, the fact
$\gamma-\frac52>0$ helps us to obtain the the desirable result
easily. Arguing similarly as before, we have
\begin{align*}
\|I_1\|_{L^1_x}&\lesssim\sum_{j\geq J-N_0}\sum_{k\leq
J-N_1}\|\triangle_k\mathcal{I}\big(
S_{j+2}{v_F}\triangle_ju\big)\|_{L^3}\|u\|_{L^\infty_T
L^6}\|u_t\|_{L^\infty_T L^2}\cr&\lesssim\sum_{j\geq
J-N_0}\sum_{k\leq
J-N_1}2^{k(\gamma-3)}2^{3k(\frac12-\frac13)}\|\triangle_k\big(
S_{j+2}{v_F}\triangle_ju\big)\|_{L^2}E_T(u)\cr&\lesssim\sum_{j\geq
J-N_0}\sum_{k\leq J-N_1}2^{k(\gamma-\frac52)}\|
S_{j+2}{v_F}\|_{L^\infty}\|\triangle_ju\|_{L^2}E_T(u).
\end{align*}
Choosing $(p,r)$ such that $\frac1p+\frac1r=\frac12$ with $2\leq
r<\infty$, the Strichartz estimate yields
\begin{align*}
\|I_1\|_{L^1_TL^1_x}&\lesssim T^{1-\frac1p}\sum_{k\leq
J-N_1}2^{k(\gamma-\frac52)}\sum_{j\geq J-N_0}2^{-j}\sum_{
j^\prime\leq
j}2^{\frac{j'}2}c_{j'}\mathcal{E}_{h,\sigma}E^{\frac32}_T(u)
\cr&\lesssim T^{1-\frac1p}2^{J(\gamma-\frac52)}\sum_{j\geq
J-N_0}2^{-\frac j2} \mathcal{E}_{h,\sigma}E^{\frac32}_T(u)
\cr&\lesssim
T^{\frac12+\frac1r}2^{-2J[s-(\frac\gamma2-\frac34+\frac1{2r})]}
\mathcal{E}^2_{s}E^{\frac32}_T(u).
\end{align*}
Combining these two cases, we have shown that
\begin{align}\label{41}
\|I_1\|_{L^1_T}\lesssim
T^{\frac12+\frac1r}2^{-2J[s-(\frac\gamma2-\frac34+\frac1{2r})]}
\mathcal{E}^2_{s}E_T(u)
\end{align}
with $\frac4{\gamma-2}\leq r<\infty$. To control
$\|I\|_{L^1_TL^1_x}$, it remains to estimate $\|I_2\|_{L^1_TL^1_x}$.
Compared with $\|I_1\|_{L^1_TL^1_x}$,  since the negative derivative
acts on the high frequency, the upper bound of
$\|I_2\|_{L^1_TL^1_x}$ is much easier to get. Here is the details:
\begin{align*}
\|I_2\|_{L^1_x}&\lesssim\sum_{j\geq J-N_0 }\sum_{k\geq
J-N_1}\|\triangle_k\mathcal{I}\big(
S_{j+2}{v_F}\triangle_ju\big)\|_{L^3}\|u\|_{L^\infty_T
L^6}\|u_t\|_{L^\infty_T L^2}\cr&\lesssim\sum_{j\geq
J-N_0}\sum_{k\geq J-N_1}2^{k(\gamma-3)}\|
S_{j+2}{v_F}\|_{L^\infty}\|\triangle_ju\|_{L^3}E_T(u).
\end{align*}
Choosing $(p,r)$ such that $\frac1p+\frac1r=\frac12$ with $2\leq
r<\infty$ again, the Strichartz estimate yields
\begin{align*}
\|I_2\|_{L^1_TL^1_x}&\lesssim T^{1-\frac1p}\sum_{k\geq
J-N_1}2^{k(\gamma-3)}\sum_{j\geq J-N_0}2^{-\frac j2}\sum_{
j^\prime\leq
j}2^{\frac{j'}2}c_{j'}\mathcal{E}_{h,\sigma}E^{\frac32}_T(u)
\cr&\lesssim
T^{\frac12+\frac1r}2^{-2J[s-(\frac\gamma2-\frac34+\frac1{2r})]}
\mathcal{E}^2_{s}E_T(u).
\end{align*}
Combining this with \eqref{41}, we obtain that
\begin{align}\label{42}
\|I\|_{L^1_TL^1_x}\lesssim
T^{\frac12+\frac1r}2^{-2J[s-(\frac\gamma2-\frac34+\frac1{2r})]}
\mathcal{E}^2_{s}E_T(u)
\end{align}
for $\frac4{\gamma-2}\leq r<\infty$.

To complete the proof the Lemma \ref{lem6}, it remains to estimate
$II$. One can proceed this as above by H\"{o}lder's inequality to
estimate
\begin{align}\label{45}
\big\|\sum_{j\geq J-N_0}2^{j(\gamma-3)}\|\triangle_jv_F
S_{j-1}{u}\|_{L^3}\big\|_{L^1_T}E_T(u).
\end{align}
Resorting to the H\"{o}lder inequality and the classical Strichartz
estimate, one can obtain that
\begin{align*}
\|II\|_{L^1_TL^1}\lesssim
T^{\frac12+\frac1r}2^{-2J[s-(\frac\gamma2-\frac34+\frac1{2r})]}
\mathcal{E}^2_{s}E_T(u).
\end{align*}
with $2\leq r\leq 6$. One also can try to improve the result by
using the precise Strichartz estimate as before, but it fails and
merely obtain that
\begin{align*}
\|II\|_{L^1_TL^1}\lesssim
T^{\frac12+\frac1r}2^{-2J[s-(\frac\gamma2-\frac34+\frac1{2r})]}
\mathcal{E}^2_{s}E_T(u).
\end{align*}
with $2\leq r\leq4$.

One can easily check that the result is worse than the desirable
result because of the restriction of $r$. Compared with the second
term in \eqref{7.0}, the negative derivative acts on the high
frequency part so that it is tempting to obtain a better result than
that of \eqref{7.0}. But $\triangle_j v_F$ is bound with $S_{j-1} u$
by the operator $\mathcal{I}$, and this structure prevents us from
using efficiently the precise Strichartz estimate. If one first
resort to the H\"{o}lder inequality, as shown in \eqref{45}, he or
she merely obtains a loss result because of the range restriction of
$r$. To go around this difficulty, we first establish a commutator
estimate through exploiting cancellation property. Now we turn to
details. Our task is to estimate
\begin{align*}\Big|\int_0^T\int_{\R^3}\sum_{j\geq J-N_0}\mathcal{I}(\triangle_jv_F
S_{j-1}{u})uu_t\mathrm{d}x\mathrm{d}t\Big|.
\end{align*}
In order to drag the $S_{j-1}u$ out of the operator $\mathcal{I}$,
we construct $u\mathcal{I}(\triangle_jv_F) S_{j-1}{u}$ and the
triangle inequality yields that
\begin{align*}\Big|\int_0^T\int_{\R^3}\sum_{j\geq J-N_0}\mathcal{I}(\triangle_jv_F
S_{j-1}{u})uu_t\mathrm{d}x\mathrm{d}t\Big|&\leq \sum_{j\geq
J-N_0}\big\|\big(\mathcal{I}(\triangle_jv_F
S_{j-1}{u})-\mathcal{I}(\triangle_jv_F)
S_{j-1}{u}\big)uu_t\big\|_{L^1_TL^1_x}\\&+\Big|\int_0^T\int_{\R^3}\sum_{j\geq
J-N_0}\mathcal{I}(\triangle_jv_F)
S_{j-1}{u}uu_t\mathrm{d}x\mathrm{d}t\Big|.
\end{align*}
We benefit from the cancellation when we deal with the first term.
Since both the Fourier transformation of $\mathcal{I}(\triangle_jv_F
S_{j-1}{u})$ and $\mathcal{I}(\triangle_jv_F) S_{j-1}{u}$ are
supported in a ring sized $2^j$, the H\"{o}lder inequality and the
Bernstein inequality lead to that
\begin{align*} \big\|\big(\mathcal{I}(\triangle_jv_F
S_{j-1}{u})-\mathcal{I}(\triangle_jv_F)
S_{j-1}{u}\big)u\big\|_{L^2_x}\leq 2^{\frac
j2}\|\mathcal{I}(\triangle_jv_F
S_{j-1}{u})-\mathcal{I}(\triangle_jv_F)
S_{j-1}{u}\|_{L^2_x}\|u\|_{L^6}.
\end{align*}
Before estimating its right hand, we recall the Coifman and Meyer
multiplier theorem. Consider an infinitely differentiable symbol
$m:\R^{nk}\mapsto\C$ so that for all $\alpha\in \N^{nk}$ and all
$\xi=(\xi_1,\xi_2,\cdots,\xi_k)\in\R^{nk}$, there is a constant
$c(\alpha)$ such that
\begin{align}\label{symbol}
|\partial_{\xi}^\alpha m(\xi)|\leq c(\alpha)(1+|\xi|)^{-|\alpha|}.
\end{align}
Define the multilinear operator $T$ by
\begin{align}
[T(f_1,\cdots,f_k)](x)=\int_{\R^{nk}}e^{ix\cdot(\xi_1+\cdots+\xi_k)}m(\xi_1,\cdots,\xi_k)\hat
f_1(\xi_1),\cdots,\hat f_k(\xi_k)d\xi_1\cdots d\xi_k,
\end{align}
or
\begin{align}
\mathcal{F}[T(f_1,\cdots,f_k)](\xi)=\int_{\xi=\xi_1+\cdots+\xi_k}m(\xi_1,\cdots,\xi_k)\hat
f_1(\xi_1),\cdots,\hat f_k(\xi_k)d\xi_1\cdots d\xi_{k-1}.
\end{align}
\begin{proposition}[\cite{C78M},Page 179.]\label{symbol1}
Suppose $p_j\in (1,\infty), j=1,\cdots k,$ are such that
$\frac1p=\frac1{p_1}+\frac1{p_2}+\cdots+\frac1{p_k}\leq1.$ Assume
$m(\xi_1,\cdots,\xi_k)$ a smooth symbol as in \eqref{symbol}. Then
there is a constant $C=C(p_i,n,k,c(\alpha))$ so that for all Schwarz
class functions $f_1,\cdots,f_k,$
\begin{align}\label{symbol2}
\|[T(f_1,\cdots,f_k)](x)\|_{L^p(\R^n)}\leq
C\|f_1\|_{L^{p_1}(\R^n)}\cdots \|f_k\|_{L^{p_k}(\R^n)}.
\end{align}
\end{proposition}
Since the operator $\mathcal{I}$ is a convolution operator with
kernel $|x|^{-\gamma}$ in $\R^3$, we can write that
$$\mathcal{F}[\mathcal{I}(\triangle_jv_F S_{j-1}{u})-\mathcal{I}(\triangle_jv_F)
S_{j-1}{u}](\xi)=\int_{\xi=\xi_1+\xi_2}\big(|\xi_1+\xi_2|^{\gamma-3}-|\xi_1|^{\gamma-3}\big)\widehat{\triangle_j
v_F}(\xi_1)\widehat{S_{j-1} u}(\xi_2)d\xi_2.$$ By the mean value
theorem, the right hand of the above formula becomes that
$$\int_{\xi=\xi_1+\xi_2}|\xi_1+\lambda\xi_2|^{\gamma-4}\frac{(\xi_1+\lambda\xi_2)\cdot \xi_2}{|\xi_1+\lambda\xi_2|}\widehat{\triangle_j
v_F}(\xi_1)\widehat{S_{j-1} u}(\xi_2)d\xi_2,$$ for a certain
$\lambda\in [0,1].$ Moreover, we rewrite it as follow:
$$\int_{\xi=\xi_1+\xi_2}m(\xi_1,\xi_2)\widehat{f_1}(\xi_1)\widehat{f_2}(\xi_2)d\xi_2,$$ with
$$m(\xi_1,\xi_2)=(\xi_1+\lambda\xi_2){|\xi_1+\lambda\xi_2|^{\gamma-5}|\xi_1|^{4-\gamma}},\quad
f_1=|\nabla|^{\gamma-4}\triangle_j v_F,\quad f_2=\nabla S_{j-1} u.$$
Observe that $|\xi_1|\geq 2^{j-1}$ and $2^{j-2}\geq |\xi_2|$, we
have that $|\xi_1+\lambda\xi_2|\sim |\xi_1|\geq 2^{J-N_0}.$ Hence,
we can check that the symbol $m(\xi_1,\xi_2)$ satisfies the estimate
\eqref{symbol}. Finally, it follows from Proposition \ref{symbol1}
that
\begin{align*} \|\mathcal{I}(\triangle_jv_F
S_{j-1}{u})-\mathcal{I}(\triangle_jv_F) S_{j-1}{u}\|_{L^2_x}\lesssim
\|f_1\|_{L^r_x}\|f_2\|_{L^{\frac{2r}{r-2}}_x}
\end{align*}
with $2<r<\infty.$ After making use of the Bernstein inequality, the
right hand can be controlled by
\begin{align*}
2^{j(\gamma-4+\frac3r)}\|\triangle_j v_F\|_{L^r_x}\|\nabla
u\|_{L^{2}_x}.
\end{align*}
Keeping in mind $j\geq J-N_0$ and recalling the definition of
$\mathcal{E}_{h,\sigma}$, the Strichartz estimate and a direct
calculation of summing in $j$ show that
\begin{align*}
T^{1-\frac1p}\sum_{j\geq J-N_0}2^{\frac
j2}2^{j(\gamma-4+\frac3r)}\|\triangle_j v_F\|_{L^p_TL^r_x}\leq
T^{1-\frac1p}\sum_{j\geq
J-N_0}2^{j(\gamma-3+\frac1r)}2^{j(\frac12-s)}\mathcal{E}_{h,s}.
\end{align*}
with $\frac1p+\frac1r=\frac12$ and $2<r<\infty.$ Choosing $r$ such
that $\max\{2,\frac1{3-\gamma}\}\leq r<\infty$, we have that
\begin{align} \sum_{j\geq
J-N_0}\big\|\big(\mathcal{I}(\triangle_jv_F
S_{j-1}{u})-\mathcal{I}(\triangle_jv_F)&
S_{j-1}{u}\big)uu_t\big\|_{L^1_TL^1_x}\cr&\lesssim
T^{\frac12+\frac1r}
2^{-2J[s-(\frac\gamma2-\frac34+\frac1{2r})]}\mathcal{E}^2_{s}E_T(u).
\end{align}
Now the rest of the paper devotes to estimate this term
\begin{align*}\Big|\int_0^T\int_{\R^3}\sum_{j\geq
J-N_0}\mathcal{I}(\triangle_jv_F)
S_{j-1}{u}uu_t\mathrm{d}x\mathrm{d}t\Big|.
\end{align*} In order to
use precise Strichartz estimate, we need to decompose this term by
Bony's para-product decomposition again,
{\begin{align*}\mathcal{I}(\triangle_jv_F) S_{j-1}{u}uu_t&=
\sum_{k}\Big\{S_{k-1}(u
S_{j-1}{u})\triangle_k\mathcal{I}(\triangle_jv_F)u_t+\triangle_{k}(u
S_{j-1}{u})S_{k+2}\mathcal{I}(\triangle_jv_F)\Big\}\cr &=II_1+II_2.
\end{align*}
After decomposing this, the term $II_1$ is similar to the second
term in the \eqref{7.0} and the negative derivative acts on the high
frequency $\triangle_j v_F$ leading to a better result than the
second term in the \eqref{7.0}. Thanks to Fourier-Plancherel formula
and H\"older inequality, we obtain
\begin{align*}
\sum_{j\geq
J-N_0}\int_0^T\int_{\R^3}II_1\mathrm{d}x\mathrm{d}t&\approx\sum_{j\geq
J-N_0}\sum_{k}\int  S_{k-1}{(u S_{j-1}u)}
\triangle_k\mathcal{I}(\triangle_jv_F)\triangle_ku_t dxdt\cr
&\approx\sum_{j\geq J-N_0}\sum_{k}\int \sum_{k'\leq k-2}
\triangle_{k'}{(u S_{j-1}u)}
\triangle_k\mathcal{I}(\triangle_jv_F)\triangle_ku_t
dxdt\cr&\lesssim\sum_{j\geq J-N_0}\sum_{k'}\int \triangle_{k'}{(u
S_{j-1}u)} \triangle_{k'}\sum_{k'\leq
k-2}(\triangle_k\mathcal{I}(\triangle_jv_F)\triangle_ku_t)
dxdt\cr&\lesssim\sum_{j\geq J-N_0}\|u S_{j-1}u\|_{L^\infty\dot
B^{\frac12}_{2,2}}\int_0^T
\big\|2^{-\frac{k'}2}\|\triangle_{k'}\sum_{k'\leq
k-2}(\triangle_k\mathcal{I}(\triangle_jv_F)\triangle_ku_t)\|_{L^2}\big\|_{\ell^2}
dt\cr &\lesssim\sum_{j\geq J-N_0}\|u\|^2_{L^\infty H^1}\int_0^T
\big\|2^{-\frac{k'}2}\|\triangle_{k'}\sum_{k'\leq
k-2}(\triangle_k\mathcal{I}(\triangle_jv_F)\triangle_ku_t)\|_{L^2}\big\|_{\ell^2}
dt
\end{align*}
On the other hand, one denotes \begin{align*}
g_{k',j}=\triangle_{k'}\sum_{k'\leq
k-2}\big(\triangle_k\mathcal{I}(\triangle_jv_F)\triangle_ku_t\big),
\end{align*}
 to estimate
$$\sum_{k'}2^{k'(-\frac12+\frac3r)}\|g_{k',j}\|_{L^1_TL^{\frac{2r}{r+2}}}.$$
Let us write that
\begin{align*}
g_{k',j}&= \sum_{k'\leq
k-2}\triangle_{k'}\bigg(\sum_{\nu\in\Lambda_{k,k'}}\triangle_{k,k'}^\nu
\mathcal{I}(\triangle_jv_F)\triangle_ku_t\bigg).
\end{align*}
As the support of the Fourier transform of a product is included in
the sum of the support of each Fourier transform, we obtain
\begin{align*}
g_{k',j}&= \sum_{k'\leq
k-2}\triangle_{k'}\bigg(\sum_{\nu\in\Lambda_{k,k'}}\triangle_{k,k'}^\nu
\mathcal{I}(\triangle_jv_F)\widetilde{\triangle_{k,k'}^\nu}
u_t\bigg).
\end{align*}
Using H\"older inequality, we get
\begin{align*}
\|g_{k',j}\|_{L^{\frac{2r}{r+2}}}&\leq \sum_{k'\leq
k-2}\sum_{\nu\in\Lambda_{k,k'}}\|\triangle_{k,k'}^\nu\mathcal{I}(\triangle_jv_F)\|_{L^r}
\|\widetilde{\triangle_{k,k'}^\nu}u_t\|_{L^{2}}\cr&\leq2^{j(\gamma-3)}
\sum_{k'\leq
k-2}\big(\sum_{\nu\in\Lambda_{k,k'}}\|\triangle_{k,k'}^\nu
v_F\|^2_{L^r}\big)^{\frac12}\big(\sum_{\nu\in\Lambda_{k,k'}}\|\triangle_{k,k'}^\nu
u_t\|^2_{L^2}\big)^{\frac12}\cr&\leq2^{j(\gamma-3)} \sum_{k'\leq
k-2}\big(\sum_{\nu\in\Lambda_{k,k'}}\|\triangle_{k,k'}^\nu
v_F\|^2_{L^r}\big)^{\frac12}\|\triangle_{k} u_t\|_{L^2}
\end{align*}
the use of quasi-orthogonality properties is made in the last
inequality.

Precise Strichartz estimate and the quasi-orthogonality properties
imply that
\begin{align*}\|g_{k',j}\|_{L^1_T(L^{\frac{2r}{r+2}})}&\leq
T^{\frac12-\frac1p} 2^{j(\gamma-3)}\sum_{k'\leq
k-2}2^{(k'-k)(\frac12-\frac1r)}2^{k(\frac32-\frac3r-\frac1p)}\cr&\,\,\,\,\,\,\,\,\,
\times \bigg(\big(\sum_{\nu\in\Lambda_{k,k'}}\|\triangle_{k,k'}^\nu
v_0\|^2_{L^2}\big)^{\frac12}+
2^{-k}\big(\sum_{\nu\in\Lambda_{k,k'}}\|\triangle_{k,k'}^\nu
v_1\|^2_{L^2}\big)^{\frac12}\bigg)\|\triangle_ku_t\|_{L^2_TL^{2}_x}\cr
\lesssim T^{\frac12-\frac1p}2^{j(\gamma-3)}& \sum_{k'\leq
k-2}2^{(k'-k)(\frac12-\frac1r)}2^{k(\frac32-\frac3r-\frac1p)}
\bigg(\|\triangle_k v_0\|_{L^2}+ 2^{-k}\|\triangle_k
v_1\|_{L^2}\bigg)\|\triangle_ku_t\|_{L^2_TL^{2}_x}
\end{align*}
with $\frac1 p+\frac1 r=\frac12$ for $2\leq r< \infty$. Therefore
\begin{align*} \sum_{k'}2^{k'(-\frac12+\frac3r)}\|g_{k',j}\|_{L^1_T(L^{\frac{2r}{r+2}})}
&\lesssim T^{\frac12-\frac1p}2^{j(\gamma-3)}\sum_{k'}\sum_{k'\leq
k-2}2^{(k'-k)\frac2r}c_k\tilde{c}_k\mathcal{E}_{h,\sigma}E^{\frac12}_T(u).
\end{align*}
A direct computation shows that
\begin{align*} \sum_{k'}2^{-\frac {k'}2}\|g_{k',j}\|_{L^1_TL^{2}}\lesssim\sum_{k'}2^{k'(-\frac12+\frac3r)}\|g_{k',j}\|_{L^1_T(L^{\frac{2r}{r+2}})}
&\lesssim
T^{\frac12-\frac1p}2^{j(\gamma-3)}\mathcal{E}_{h,\sigma}E^{\frac12}_T(u).
\end{align*}
Hence, we have that
\begin{align}\label{46}
\Big|\sum_{j\geq
J-N_0}\int_0^T\int_{\R^3}II_1\mathrm{d}x\mathrm{d}t\Big|\lesssim2^{-2J[s-(\frac\gamma2-\frac34+\frac1{2r})]}T^{\frac12+\frac1r}\mathcal{E}^2_{s}E_T(u)
\end{align}
with $\frac4{\gamma-2}\leq r<\infty$. Finally, we conclude this
section by giving the estimate of $II_2$.
\begin{align}\label{47}\Big|\sum_{j\geq
J-N_0}\int_0^T\int_{\R^3}II_2\mathrm{d}x\mathrm{d}t\Big|&\lesssim
T^{\frac12}\sum_{j\geq J-N_0}\sum_{k}\|\triangle_{k}(u
S_{j-1}{u})S_{k+1}\mathcal{I}(\triangle_jv_F)\|_{L_T^2L^2}\|u_t\|_{L^\infty_T
L^2}\cr\lesssim T^{1-\frac1p}\sum_{j\geq
J-N_0}2^{j(\gamma-3)}&\sum_{ k}\sum_{ k'\leq k}\|\triangle_{k}(u
S_{j-1}{u})\|_{L_T^\infty
L^2}2^{k'\frac3r}\|\triangle_{k'}\triangle_jv_F\|_{L_T^pL^r}E^{\frac12}_T(u)\cr\lesssim
T^{1-\frac1p}\sum_{j\geq J-N_0}2^{j(\gamma-3)}&\sum_{ k}\sum_{k'\leq
k}\|\triangle_{k}(u S_{j-1}{u})\|_{L_T^\infty
L^2}2^{\frac{k'}2}c_{k'}\mathcal{E}_{h,\sigma}E^{\frac12}_T(u)\cr\lesssim
T^{1-\frac1p}\sum_{j\geq J-N_0}2^{j(\gamma-3)}&\sum_{
k'}c_{k'}\sum_{k'\leq k}2^{\frac k 2}\|\triangle_{k}(u
S_{j-1}{u})\|_{L_T^\infty
L^2}2^{(k'-k)\frac12}\mathcal{E}_{h,\sigma}E^{\frac12}_T(u)\cr
\lesssim T^{1-\frac1p}\sum_{j\geq
J-N_0}2^{j(\gamma-3)}&\|c_{k'}\|_{\ell^2(\Z)}\|2^{\frac k
2}\|\triangle_{k}(u S_{j-1}{u})\|_{L_T^\infty
L^2}\|_{\ell^2(\Z)}\|2^{-\frac{k}2}\|_{\ell^2(\N)}\mathcal{E}_{h,\sigma}E^{\frac12}_T(u)\cr\lesssim
T^{1-\frac1p}2^{J(\gamma-3)}E^{\frac32}_T(u)&\mathcal{E}_{h,\sigma}\lesssim2^{-2J[s-(\frac\gamma2-\frac34+\frac1{2r})]}T^{\frac12+\frac1r}\mathcal{E}^2_{s}E_T(u).
\end{align}
Collecting \eqref{46} and \eqref{47}, we have been proved that
\begin{align}\label{48}
\Big|\int_0^T\int_{\R^3}\sum_{j\geq
J-N_0}\mathcal{I}(\triangle_jv_F)
S_{j-1}{u}uu_t\mathrm{d}x\mathrm{d}t\Big|\lesssim
T^{\frac12+\frac1r}2^{-2J[s-(\frac\gamma2-\frac34+\frac1{2r})]}
\mathcal{E}^2_{s}E_T(u),
\end{align}
with $\frac4{\gamma-2}\leq r<\infty$. Finally, we complete the
proof of \eqref{34} by \eqref{42} and \eqref{48}, hence it ends
the proof of Lemma \ref{lem6}.

\vskip0.5cm {\bf Acknowledgements:}\quad The authors are grateful
to  Prof. J.Chemin for sending his lecture to us.   The authors
were partly supported by the NSF of China, No.10725102.

\begin{center}

\end{center}
\end{document}